%% file: main.tex
\documentclass[preprint,11pt,authoryear,nonatbib]{elsarticle}
\usepackage[margin=1in]{geometry}

\usepackage{amssymb}
\usepackage{placeins}
\usepackage{algorithm}%
\usepackage{algorithmicx}%
\usepackage{algpseudocode}%

\usepackage{graphicx}

\usepackage{amsmath}
\usepackage{bm}
\usepackage{tabularx}
\usepackage{booktabs}
\usepackage{array}
\usepackage{comment}
\usepackage{xcolor}
\usepackage{amsthm}
\newtheorem{definition}{Definition}
\newtheorem{remark}{Remark}
\usepackage{setspace}
\usepackage[caption=false,font=footnotesize]{subfig}
\input{mycommands}

\journal{European Journal of Operational Research}

\usepackage{csquotes}
\usepackage[
    backend=biber,
    style=apa,
    natbib=true
]{biblatex}

\begin{document}
\onehalfspacing
\begin{frontmatter}



\title{Computing Equilibria in Simulation-Based Insurance Markets with Discontinuous Demand} 


\author{Yunsoo Ha, Linda Nozick} 

\affiliation{organization={School of Civil and Environmental Engineering, Cornell University},
            addressline={Hollister Hall}, 
            city={Ithaca},
            postcode={14853}, 
            state={NY},
            country={United States}}

\begin{abstract}
We study equilibrium computation in simulation-based games motivated by competition in natural-hazard insurance markets. The motivating application considers a multi-region market in which each insurer selects regional pricing and reinsurance decisions to maximize expected financial performance subject to multi-period insolvency constraints. Competitors' decisions affect both profitability and insolvency feasibility through the insured portfolio. Consequently, the strategic interaction defines a generalized Nash game with discontinuous player objectives induced by discrete household purchase decisions and potentially nonconvex feasible sets generated by the insolvency constraints. Discontinuity, nonconvex feasibility, and simulation-defined objectives and constraints complicate the direct use of standard optimality-based equilibrium formulations. Moreover, repeated evaluation of the multi-period simulator makes equilibrium search computationally expensive. To address these computational challenges, we develop an evaluation-efficient framework for computing approximate local equilibria, which combines a response optimizer tailored to the distinct structures of pricing and reinsurance with damped better-response dynamics to mitigate cycling. In a case study of the North Carolina hurricane insurance market, the framework identifies multiple approximate local equilibrium candidates both with and without a household affordability cap and finds no profitable sampled unilateral deviations within the tested neighborhoods under the original simulation model. A representative equilibrium candidate is computed in approximately 5--10 minutes on a single GPU, making repeated equilibrium analyses computationally practical.
\end{abstract}



\begin{keyword}
Simulation-based games \sep Equilibrium computation\sep Derivative-free optimization\sep Discontinuous optimization\sep Natural-hazard insurance


\end{keyword}

\end{frontmatter}

\section{Introduction}
\label{sec:introduction}

Natural-hazard insurance provides households with financial protection against disaster losses, while requiring insurers to manage potentially large and uncertain claims. To manage the financial risk associated with disaster claims, insurers must coordinate pricing and reinsurance decisions. The premiums set by an insurer affect which households purchase coverage, determining both premium income and the insurer's exposure to disaster losses. Meanwhile, reinsurance transfers part of the loss exposure to a reinsurer at a cost that depends on the coverage purchased and the insured risks. Consequently, the profitability of a pricing decision and the insurer's ability to satisfy solvency requirements depend on both the households attracted and the reinsurance coverage selected. At the same time, competition introduces strategic interdependence among insurers because households compare premiums across insurers when choosing coverage. Hence, understanding the resulting market outcomes requires an equilibrium analysis of strategic interactions among insurers in pricing and reinsurance.

Established operations research methods for computing Nash and generalized Nash equilibria offer a natural starting point for analyzing competition among insurers. In particular, complementarity formulations characterize equilibria by combining the optimality conditions of the individual players’ problems into a system of complementarity conditions, as illustrated by the Cournot--Nash computations in~\citet{kolstad1991computing}. Generalized Nash equilibrium methods further accommodate strategy-dependent feasible sets, in which the feasible decisions of one player depend on the actions of other players~\citep{facchinei2010generalized,migot2020parametrized}. However, the applicability of optimality-based equilibrium methods depends on the structural properties of the underlying player problems. In the natural-hazard insurance market considered in this study, discrete household participation and simulation-based solvency constraints produce discontinuous objective functions and potentially nonconvex feasible strategy sets, complicating the direct use of standard optimality-based formulations.

To address these difficulties, we develop a computational framework for a multi-region catastrophe insurance market in which insurers jointly choose regional pricing parameters and excess-of-loss reinsurance contracts. The model combines household purchase decisions based on utility comparisons and, when imposed, affordability thresholds with a multi-period simulation of insurer capital. Competitors' decisions affect both the objective value and the feasible strategy set through the portfolio available to each insurer, giving rise to a generalized Nash equilibrium problem. The response procedure employs complementary search strategies tailored to the distinct computational characteristics of pricing and reinsurance decisions, while evaluating candidate improvements using the original simulation. More broadly, the proposed framework can accommodate \textit{simulation-based games with discontinuous objectives and nonconvex feasible sets}, with the specific search strategies adapted to the structure of the decision variables. Our focus is on computing approximate local equilibrium candidates efficiently enough to support repeated analyses of market conditions and policy choices.

\subsection{Literature Review}
\label{subsec:literature_positioning}

The computational task has two connected levels: identifying
a market equilibrium and optimizing an individual insurer's
response to its competitors. We position the proposed framework at both levels. Table~\ref{tab:equilibrium_literature} summarizes
representative approaches to equilibrium computation and identifies barriers to direct application in the present model.
Table~\ref{tab:response_optimization_literature} reviews optimization methods relevant to the insurer response problem and identifies the additional treatment required by its computational structure.

\subsubsection{Equilibrium Computation and Insurance-Market Models}
\label{subsubsec:equilibrium_literature}

\medskip
\noindent
\textbf{Equilibrium computation with nonconvex player problems.}
Local optimality conditions provide one route to equilibrium characterization when player problems are nonconvex. \citet{pang2011nonconvex} formulate quasi-Nash equilibria through aggregated first-order conditions and establish second-order sufficient conditions guaranteeing that a quasi-Nash equilibrium is a local Nash equilibrium. However, payoff discontinuities can make first-order conditions unavailable or insufficient for equilibrium characterization. \citet{feinstein2024characterizing} provide a global characterization of the Nash equilibrium set through the Pareto-optimal solutions of a vector optimization problem with a nonconvex ordering cone. Related vector optimization characterizations provide necessary and sufficient conditions for generalized Nash equilibria, and the characterization results require neither smoothness nor convexity. Nevertheless, computational implementation requires solving nonconvex vector optimization problems, leaving practical tractability dependent on additional problem structure. Stronger algorithmic conclusions are available under more restrictive structural assumptions. \citet{carvalho2026cut} develop Cut-and-Play for nonconvex games with separable payoffs and polyhedrally representable feasible sets. The algorithm iteratively refines convex outer approximations using cutting planes and branching to compute mixed-strategy Nash equilibria or certify nonexistence. Although continuity and convexity are not required, applicability depends on payoff separability and polyhedral representability.

The preceding literature shows that reformulation-based equilibrium computation is most effective when player response problems possess analytically exploitable structure. However, the insurer response problem considered here does not readily admit such a representation.
Pricing and reinsurance are chosen jointly, price changes produce discrete shifts in household participation and the insured portfolio, and competitors' prices alter the demand underlying insolvency feasibility. 
Consequently, an insurer's objective value and satisfaction of the insolvency constraint can be evaluated only through the multi-period capital simulation.
Accordingly, equilibrium computation must be built around
simulation-based optimization of individual insurer responses. 
\medskip

\noindent
\textbf{Simulation-based and empirical game analysis.}
Simulation-based game analysis uses simulator output to support strategy exploration, payoff estimation, and equilibrium computation~\citep{wellman2025empirical}. For simulation-based games with finite-dimensional continuous strategy spaces, \citet{vorobeychik2008stochastic} develop stochastic-search procedures that directly approximate best responses and Nash equilibria. However, direct search over continuous strategy profiles can require a large number of simulator evaluations. To manage this cost, restricted empirical-game methods compute equilibria over finite sets of evaluated strategies and selectively expand the restricted games through searches for profitable deviations. \citet{jordan2010strategy} study strategy-exploration policies for such iterative expansion, including the selection of improving deviations and best responses to equilibria of the current restricted game. 
Another way to reduce the evaluation burden is to predict payoffs at unevaluated profiles from accumulated simulator observations. \citet{picheny2019bayesian} develop a Gaussian-process approach for equilibrium computation, using predictive means and uncertainties to guide the selection of additional simulations. Sampling criteria based on equilibrium probability or uncertainty reduction concentrate the evaluation budget on profiles that are informative for identifying Nash equilibria. However, reducing the evaluation burden does not address the additional challenges posed by payoff discontinuities.

\medskip

\noindent
\textbf{Discontinuous games.}
Discontinuous games have received considerable theoretical attention. \citet{reny1999existence} provide sufficient conditions for the existence of pure-strategy Nash equilibria in discontinuous games. Subsequent research develops additional existence results for discontinuous games with complete and incomplete information; see, for example, \citet{carbonell2018existence} and the review by \citet{reny2020nash}. Taken together, the existence results show that payoff discontinuities do not rule out equilibrium existence when appropriate structural conditions hold. However, conditions ensuring equilibrium existence in discontinuous games can be difficult to verify and do not directly yield tractable computational procedures~\citep{nessah2008existence}. Accordingly, we focus on developing a computational procedure for finding local equilibria in the simulation-based insurance game. The main computational difficulty arises from discontinuities generated endogenously by household participation decisions: crossing a household-specific threshold changes the insured portfolio and, consequently, simulated losses, premium revenues, and solvency. The consequences of participation-induced discontinuities depend on the modeling of household demand and insurer risk. Accordingly, we next review demand modeling and related research on insurance-market equilibria.
\medskip

\noindent
\textbf{Demand modeling and insurance-market equilibria.}
In price competition, \citet{liu2025bayesian} study Bayesian Nash equilibria under multinomial logit demand, which yields smooth purchase probabilities as functions of firms' prices. This smooth demand structure enables a stochastic variational-inequality formulation and a finite-dimensional approximation, in contrast to the discrete household purchase indicators in the setting considered here. Beyond the pricing considerations in standard price competition, an insurer's response problem must also account for the risk associated with the business attracted by its pricing decision. \citet{dutang2013competition} incorporate this connection through premium competition under solvency constraints. Strategic interactions also depend on the timing of insurers' decisions, as insurers may choose their actions simultaneously or respond sequentially to competitors. Reflecting this distinction, \citet{BOONEN2026972} examine both Nash and Stackelberg equilibria in a single-period stochastic insurance duopoly.

For catastrophe insurance across multiple regions, pricing also determines the geographical composition of the portfolio and its associated reinsurance requirements. This connection is central to \citet{liu2025computing}, the closest application antecedent, who study multi-region hurricane-insurance competition with joint pricing and reinsurance decisions. We build on that market setting by incorporating multi-period capital dynamics and controlling the expected proportion of periods in which insolvency occurs. The resulting feasibility assessment depends on the evolution of capital along simulated catastrophe paths. Our computational contribution is to exploit the pricing--reinsurance structure within these response problems and use the resulting responses to search for symmetric approximate local equilibria in pure strategies. 

\begin{table}[!htbp]
\centering
\caption{Equilibrium approaches and barriers to direct application.}
\label{tab:equilibrium_literature}

\begingroup
\setstretch{1}
\small

\setlength{\tabcolsep}{6pt}
\renewcommand{\arraystretch}{1.2}
\setlength{\extrarowheight}{2pt}

\setlength{\aboverulesep}{0.5ex}
\setlength{\belowrulesep}{0.5ex}

\renewcommand{\tabularxcolumn}[1]{m{#1}}

\begin{tabularx}{\linewidth}{@{}
    >{\raggedright\arraybackslash}m{0.3\linewidth}
    >{\raggedright\arraybackslash}m{0.23\linewidth}
    >{\raggedright\arraybackslash}X@{}}

\toprule
\textbf{Reference}
& \textbf{Approach}
& \textbf{Barrier in the present model} \\
\midrule

\citet{kolstad1991computing}
& Cournot--Nash; complementarity
& Smooth player optimality conditions do not cover
participation-induced payoff jumps. \\
\midrule[0.25pt]

\citet{facchinei2010generalized};\newline
\citet{migot2020parametrized}
& GNE methods; parametrized VIs
& A tractable VI formulation of the discontinuous,
simulation-constrained responses is still needed. \\
\midrule[0.25pt]

\citet{pang2011nonconvex}
& Quasi-NE; first-order conditions
& Derivatives within a participation region do not
assess gains from crossing its boundaries. \\
\midrule[0.25pt]

\citet{feinstein2024characterizing}
& NE and GNE; vector optimization
& The reformulation still requires solving
discontinuous, simulation-defined optimization problems. \\
\midrule[0.25pt]

\citet{carvalho2026cut}
& Separable-payoff mixed NE; convexification
& Rival-dependent solvency constraints fall outside
the independent-feasible-set formulation. \\
\midrule[0.25pt]

\citet{vorobeychik2008stochastic};\newline
\citet{wellman2025empirical}
& Simulation-based search; empirical game analysis
& Requires a pricing--reinsurance response solver
with solvency constraints. \\
\midrule[0.25pt]

\citet{picheny2019bayesian}
& GP-based equilibrium search
& Requires treatment of payoff jumps and
rival-dependent feasibility. \\
\midrule[0.25pt]

\citet{liu2025bayesian}
& Logit demand; stochastic VI approximation
& The smooth-demand formulation does not retain
household participation jumps. \\
\midrule[0.25pt]

\citet{dutang2013competition};\newline
\citet{BOONEN2026972}
& Insurance pricing; solvency and strategic timing
& Extension is needed to joint pricing--reinsurance
responses with simulated insolvency-frequency constraints. \\
\midrule[0.25pt]

\citet{liu2025computing}
& Multi-region competition; pricing and reinsurance
& Computationally costly for repeated equilibrium analyses;
requires extension to multi-period capital dynamics.\\

\bottomrule
\end{tabularx}

\par\medskip
\begin{minipage}{\linewidth}
\footnotesize
\raggedright
\textit{Notes:}
NE = Nash equilibrium;
GNE = generalized Nash equilibrium;
GP = Gaussian process;
VI = variational inequality.
Barriers refer to unmet application conditions or additional
modeling and computational requirements.
\end{minipage}

\endgroup
\end{table}

\subsubsection{Optimization of Individual Insurers' Responses}
\label{subsubsec:response_optimization_literature}

With competitors' decisions fixed, each insurer faces a pricing and reinsurance problem with a discontinuous objective and a nonconvex feasible region. Motivated by these features, we begin by reviewing methods for nonsmooth nonconvex optimization and discuss their applicability to the insurer's response problem. 

\medskip
\noindent
\textbf{Nonsmooth and derivative-free local methods.}
Nonsmooth nonconvex optimization methods seek useful search information even where ordinary gradients do not exist. Subgradient methods, bundle methods, gradient sampling, and smoothing techniques address this difficulty through generalized derivatives, gradients at nearby differentiable points, or approximations of the objective~\citep{clarke1990optimization,burke2005robust,curtis2025practical}.
When even such derivative information is unavailable or prohibitively
expensive to obtain, derivative-free methods instead construct local
models or explore trial points using only function evaluations~\citep{larson2019derivative,liuzzi2019trust}.
However, the difficulty in our response problem goes beyond the absence of ordinary derivatives. Threshold-based household participation creates actual jumps in the objective, violating the local Lipschitz continuity required by many nonsmooth optimization methods, including some derivative-free approaches. In particular, finite-difference estimates that cross a participation boundary can be highly sensitive to the perturbation size. Hence, discontinuities complicate both the construction of reliable search directions and the application of standard stationarity conditions, motivating methods specifically designed for discontinuous
optimization~\citep{vicente2012discontinuous,cui2023minimization}.

\medskip
\noindent
\textbf{Smoothing and structured discontinuous optimization.}
One possible approach to discontinuous optimization is to replace the original problem with a sequence of smooth approximations, allowing gradient information to guide the search. Under suitable conditions, limits of stationary points of the approximating problems satisfy a generalized stationarity condition for the original problem~\citep{ermoliev1995minimization}.
In our setting, smoothing household participation provides local search information, but its effect extends to both profitability and risk. Replacing discrete purchase indicators with smooth approximations changes the insured portfolio and its associated cash flows, potentially altering both response incentives and the equilibrium structure. Consequently, feasibility under the smoothed model need not imply feasibility under the original insolvency constraint.
Moreover, decreasing the smoothing parameter sharpens the transition near participation thresholds, which can make gradient-based updates increasingly sensitive in these regions. 

Along a related line of work, \citet{cui2023minimization} introduce \textit{pseudo B-stationarity} for discontinuous piecewise optimization.
Pseudo B-stationarity requires B-stationarity of an auxiliary problem that fixes each indicator value according to the sign of the corresponding step-function argument at the point under consideration.
However, in the present setting, pseudo B-stationarity alone provides limited information about profitable unilateral deviations because the associated auxiliary problem considers only local perturbations that preserve the household participation pattern implied by the prospective solution. From an economic perspective, a central question is whether an insurer can increase profit by changing its price and attracting or losing households. Consequently, price changes that alter household participation fall outside the auxiliary problem, meaning that a pseudo B-stationary point may still admit a profitable, insolvency-feasible pricing deviation.

\medskip
\noindent
\textbf{Direct search for discontinuous objectives.}
An economic assessment of equilibrium requires checking whether an insurer can increase its payoff through feasible unilateral price changes, including changes that alter household participation.
Direct-search methods offer a natural approach to examining such deviations by comparing objective values at trial decisions without requiring derivatives~\citep{kolda2003optimization}.
Evaluating candidate prices under the original purchase rules allows profit gains from attracting or losing households to be assessed alongside insolvency feasibility.
However, convergence analysis for discontinuous objectives still requires careful treatment of behavior near the jumps. \citet{vicente2012discontinuous} establish generalized stationarity results for direct search under suitable conditions on discontinuous objectives. A subsequent counterexample shows that the required function-value convergence is not guaranteed under previously proposed conditions, motivating an additional revealing poll step that samples neighboring continuity regions, with related convergence properties recovered under further assumptions~\citep{audet2024counterexample}.
Beyond convergence considerations, relying exclusively on the original simulator can make response optimization computationally costly. In the insurance model, each candidate response requires simulator-based assessments of both profitability and insolvency feasibility. Hence, practical equilibrium computation requires a response optimizer that uses simulator evaluations efficiently.

\begin{table}[!htbp]
\centering
\caption{Optimization approaches and barriers in the present model.}
\label{tab:response_optimization_literature}

\begingroup
\setstretch{1}
\small

\setlength{\tabcolsep}{6pt}
\renewcommand{\arraystretch}{1.2}
\setlength{\extrarowheight}{2pt}

\setlength{\aboverulesep}{0.5ex}
\setlength{\belowrulesep}{0.5ex}

\renewcommand{\tabularxcolumn}[1]{m{#1}}

\begin{tabularx}{\linewidth}{@{}
    >{\raggedright\arraybackslash}m{0.27\linewidth}
    >{\raggedright\arraybackslash}m{0.23\linewidth}
    >{\raggedright\arraybackslash}X@{}}

\toprule
\textbf{Reference}
& \textbf{Approach}
& \textbf{Barrier in the present model} \\
\midrule

\citet{burke2005robust}
& Gradient sampling
& Participation jumps violate local Lipschitz continuity. \\
\midrule[0.25pt]

\citet{liuzzi2019trust}
& Nonsmooth derivative-free trust regions
& Locally Lipschitz analysis does not cover
participation jumps. \\
\midrule[0.25pt]

\citet{ermoliev1995minimization}
& Smooth approximations
& Smoothed feasibility need not imply
original-model feasibility. \\
\midrule[0.25pt]

\citet{cui2023minimization}
& Pseudo B-stationarity
& Does not rule out profitable deviations
that change participation. \\
\midrule[0.25pt]

\citet{kolda2003optimization}
& Direct search
& Discontinuities and insolvency feasibility
require additional treatment. \\
\midrule[0.25pt]

\citet{vicente2012discontinuous};\newline
\citet{audet2024counterexample}
& Discontinuous direct search; revealing poll
& Does not address nonconvex, simulation-defined
feasibility or repeated evaluation-efficient response search. \\

\bottomrule
\end{tabularx}

\par\medskip
\begin{minipage}{\linewidth}
\footnotesize
\raggedright
\textit{Notes:}
Barriers indicate analytical limitations or additional
application requirements.
\end{minipage}

\endgroup
\end{table}

Taken together, the two strands of literature motivate an integrated and computationally practical framework for equilibrium computation and response optimization in simulation-based games with discontinuous player problems, which forms the basis for the contributions summarized below.

\subsection{Contributions}
\label{subsec:contributions}

The paper makes the following contributions.

\begin{enumerate}

    \item
    \textbf{A simulation-based insurance equilibrium model with
    multi-period solvency constraints.}
    Extending the multi-region pricing and reinsurance setting of
    \citet{liu2025computing}, we model multi-period capital evolution
    and constrain the expected proportion of insolvent periods.
    Household participation thresholds generate payoff discontinuities,
    while competitors' decisions affect profitability and insolvency
    feasibility through the insured portfolio. The formulation yields
    a generalized Nash equilibrium problem with discontinuous objectives
    and simulation-defined, nonconvex feasible sets.

    \item
    \textbf{An evaluation-efficient framework for approximate local equilibrium search in discontinuous simulation-based games.}
    We develop a response-based framework for simulation-based games in which repeated player optimization is complicated by discontinuous payoffs, nonconvex feasibility, and expensive simulator evaluations.
    In the insurance application, the framework combines smoothing-guided pricing search, derivative-free local refinement, and tailored reinsurance optimization to reduce the simulator evaluations required for computing approximate better responses. The resulting responses are embedded within a damped symmetric equilibrium iteration, while candidate solutions are assessed under the original discontinuous simulation model.
    
    \item
    \textbf{Deviation-based assessment of equilibrium quality.}
   We assess candidate solutions by sampling unilateral pricing deviations within prescribed radii while holding competitors' decisions fixed.
    Reinsurance is approximately reoptimized for each sampled price, and payoff improvement and insolvency feasibility are evaluated under the original simulation model. Hence, the assessment captures deviations
    that change household participation and provides numerical evidence of
    local equilibrium quality within the scope of the sampled search.

    \item
    \textbf{Computational evidence and market insights.}
    We study four homogeneous insurers operating across eight regions
    in North Carolina over a thirty-year horizon. The experiments identify
    multiple approximate local equilibrium candidates both with and without
    a 5\% household affordability cap, with no profitable deviations found
    in the reported deviation tests. A single equilibrium candidate can be
    computed in approximately 5--10 minutes on a single GPU. A representative
    computation requires approximately 9,900 pricing-search objective
    evaluations, compared with at least 190,000 evaluations reported by
    \citet{liu2025computing}. The resulting candidates exhibit distinct
    combinations of premiums, reinsurance protection, and insurer financial
    performance, highlighting both the computational scalability of the
    framework and the importance of accounting for multiple possible market
    outcomes.

\end{enumerate}

The remainder of the paper presents the insurance-market model and equilibrium criteria, develops the computational framework, reports the numerical experiments and deviation tests, and concludes with directions for further research.

\section{Market Model and Equilibrium Problem}
\label{sec:model}

We consider a multi-region catastrophe insurance market in which
insurers jointly select regional pricing parameters and reinsurance
contracts. The model builds on \citet{liu2025computing} by incorporating
multi-period capital evolution and an insolvency-frequency constraint.
We describe the household demand and financial dynamics needed to
formulate each insurer's response problem, define the equilibrium
criteria, and identify the resulting computational structure.
Additional modeling details are provided in \ref{appdx:s1} and \ref{appdx:s2}.

\subsection{Insurers' Decisions and Household Demand}
\label{subsec:decisions_demand}

Consider $N$ homogeneous insurers operating across $d$ regions.
Insurer $i$ chooses a regional pricing vector
$\boldsymbol{\lambda}_i\in\Lambda\subset\mathbb{R}^d$ and an
excess-of-loss reinsurance contract parameterized by its attachment
point $A_i$ and exhaustion point $M_i$. The pricing set $\Lambda$
is a prescribed box, while the admissible reinsurance set
$\mathcal R\subset\mathbb{R}^2$ imposes parameter bounds and
$0\leq A_i<M_i$. Define
\(
    \BFx_i:=(\boldsymbol{\lambda}_i,A_i,M_i)
    \in\mathcal X:=\Lambda\times\mathcal R, 
\)
\(
    \BFx:=(\BFx_1,\ldots,\BFx_N),
\)
and let $\BFx_{-i}$ collect the decisions of all insurers other
than $i$. Insurers are homogeneous in their demand, cost, and risk-evaluation models, but different pricing decisions can result in different insured portfolios.

Household demand follows individual utility comparisons between
purchasing insurance and remaining uninsured. Let
$P_j(\boldsymbol{\lambda})$ denote the annual premium relevant
to household $j$'s purchase decision, and let
$\mathcal A_j\subseteq\Lambda$ denote the pricing vectors under
which purchasing insurance is preferred according to the utility
comparison. When an affordability cap is imposed, define
\begin{equation}
    \mathcal B_j
    :=
    \left\{
        \boldsymbol{\lambda}\in\Lambda:
        P_j(\boldsymbol{\lambda})\leq\alpha H_j
    \right\},
    \label{eq:affordability}
\end{equation}
where $H_j$ is the household's property value and
$\alpha\in(0,1)$ is the affordability ratio. Without the cap,
set $\mathcal B_j=\Lambda$. Then the purchase indicator becomes
\(
    b_j(\boldsymbol{\lambda})
    :=
    \mathbf 1
    \left\{
        \boldsymbol{\lambda}\in
        \mathcal A_j\cap\mathcal B_j
    \right\}. 
\)
Coverage-type indices are suppressed here; the wind and flood purchase rules and the application of the household affordability condition are specified in~\ref{appdx:s1}. 

Aggregating the household purchase indicators produces total market
demand across regions and coverage types. Because the household
population is finite, total demand changes discretely when prices
cross household-specific purchase thresholds. During insurer $i$'s
response optimization, the demand currently served by competing
insurers is held fixed. For each candidate pricing vector
$\boldsymbol{\lambda}_i$, total market demand is evaluated from the
household purchase rules, and insurer $i$ serves the residual demand
remaining after subtracting the demand served by its competitors.
The resulting insurer-specific demand determines premium income and
insured loss exposure and therefore affects insolvency risk.
Reinsurance does not enter the household purchase rules but adjusts
the portion of the resulting loss exposure retained by the insurer.

\subsection{Capital Dynamics and Insolvency Constraints}
\label{subsec:capital_insolvency}

Catastrophe uncertainty is represented by a sequence of events
$\xi$ over a $T$-year horizon. Following the hazard setting of
\citet{apivatanagul2011long}, the numerical model uses 97 hurricane
event types to generate scenario paths, each of which can contain
multiple events within a year. Let $\mathcal S$ denote the finite
set of paths used in the computational model, with probabilities
$p_s\geq0$ satisfying $\sum_{s\in\mathcal S}p_s=1$. The same paths
and probabilities are used when comparing strategies.

For strategy profile $\BFx$, let $C_{i,y}^s(\BFx)$ denote insurer
$i$'s capital at the end of year $y$ in scenario $s$, after any
external capital injection. Before an injection, capital evolves as
\begin{equation}
\begin{aligned}
\widetilde C_{i,y}^s(\BFx)
={}&
C_{i,y-1}^s(\BFx)
+\beta^{y-1}
\bigl[
\Pi_i(\BFx)
+D_{i,y}^s(\BFx)
-L_{i,y}^s(\BFx)
+R_{i,y}^s(\BFx)
\bigr],
\end{aligned}
\label{eq}
\end{equation}
where $\beta\in(0,1]$ is the discount factor. Here,
$\Pi_i$ is annual premium income net of operating costs,
$D_{i,y}^s$ is the deductible amount retained by policyholders,
$L_{i,y}^s$ is the insured portfolio loss before deductibles,
and $R_{i,y}^s$ is the net reinsurance cash flow, incorporating
recoveries from the excess-of-loss layer $[A_i,M_i]$ and the
associated reinsurance costs. Thus,
$L_{i,y}^s-D_{i,y}^s$ is the insurer's claim payment before
reinsurance. The construction of these cash-flow components,
including the pricing of reinsurance and the aggregation of
event-level recoveries, is provided in \ref{appdx:s2}.

An insolvency event occurs when
$\widetilde C_{i,y}^s(\BFx)\leq0$. The incumbent insurer then exits,
and a new insurer enters and continues under the same pricing and
reinsurance strategy. The entry is represented by an external capital
injection equal to three times annual premium income net of operating
costs:
\begin{align*}
    I_{i,y}^s(\BFx)
    &:=
    3\bigl(\Pi_i(\BFx)\bigr)_+
    \mathbf 1
    \left\{
        \widetilde C_{i,y}^s(\BFx)\leq0
    \right\},\\
    C_{i,y}^s(\BFx)
    &:=
    \widetilde C_{i,y}^s(\BFx)
    +I_{i,y}^s(\BFx),
\end{align*}
where $(a)_+:=\max\{a,0\}$. The simulation subsequently continues
with the new insurer.

For scenario $s$, insurer $i$'s financial performance is defined as
terminal capital net of the external capital supplied to replacement
insurers,
$F_i(\BFx,s):=C_{i,T}^s(\BFx)-\sum_{y=1}^{T}I_{i,y}^s(\BFx)$.
Subtracting the injections ensures that externally supplied capital
is not counted as performance generated by the strategy. As a result, insurer $i$'s expected payoff becomes
\[
f_i(\BFx_i,\BFx_{-i})
:=\mathbb E_{\xi}\!\left[F_i(\BFx,\xi)\right]
=\sum_{s\in\mathcal S}p_s\left(C_{i,T}^s(\BFx)-\sum_{y=1}^{T}I_{i,y}^s(\BFx)\right).
\]
Solvency is controlled through the expected proportion of years in
which insolvency occurs,
\[
g_i(\BFx_i,\BFx_{-i})
:=
\mathbb E_{\xi}
\left[
    \frac{1}{T}
    \sum_{y=1}^{T}
    \mathbf 1
    \left\{
        \widetilde C_{i,y}(\BFx,\xi)\leq0
    \right\}
\right]
=
\frac{1}{T}
\sum_{s\in\mathcal S}p_s
\sum_{y=1}^{T}
\mathbf 1
\left\{
    \widetilde C_{i,y}^s(\BFx)\leq0
\right\}.
\]
An admissible response satisfies
$g_i(\BFx_i,\BFx_{-i})\leq\varepsilon$, where $\varepsilon$ is the
prescribed insolvency tolerance. The constraint controls
the expected proportion of years in which insolvency occurs, rather
than the probability of at least one insolvency over the entire
horizon. Hence, evaluating feasibility requires simulation of the
complete multi-period capital path.

\subsection{Player Problems and Equilibrium Criteria}
\label{sec:intro-LNE}

The expected financial objective and insolvency-frequency constraint
define insurer $i$'s response problem for fixed competing strategies
$\BFx_{-i}$ as follows.
\begin{equation}
\begin{aligned}
    \max_{\BFx_i\in\mathcal X}\quad
    & f_i(\BFx_i,\BFx_{-i})\\
    \text{s.t.}\quad
    & g_i(\BFx_i,\BFx_{-i})\leq\varepsilon.
\end{aligned}
\label{eq:br}
\end{equation}
Define its feasible strategy set as
\(
    \widetilde{\mathcal X}_i(\BFx_{-i})
    :=
    \left\{
        \BFx_i\in\mathcal X:
        g_i(\BFx_i,\BFx_{-i})\leq\varepsilon
    \right\}.
\)
Competitors' decisions affect both the objective and the feasible set through the insured portfolio. Hence, the market gives rise to a generalized Nash equilibrium problem.

A feasible strategy profile $\BFx^*$ is a generalized Nash equilibrium
if $\BFx_i^*$ is a global maximizer of insurer $i$'s response problem
for every $i$. Computing a generalized Nash equilibrium requires finding
a profile at which all insurers globally optimize their responses,
while verifying a candidate requires ruling out every feasible and
profitable unilateral deviation. Both tasks are computationally
prohibitive for the simulation-based, discontinuous, and nonconvex
response problems considered here. We therefore formulate and
numerically assess local equilibrium criteria.

\begin{definition}[Local generalized Nash equilibrium]
\label{def:local-nash}
A strategy profile
$\BFx^*=(\BFx_1^*,\ldots,\BFx_N^*)$ is a local generalized Nash
equilibrium if, for every insurer $i$,
$\BFx_i^*\in\widetilde{\mathcal X}_i(\BFx_{-i}^*)$ and there exists
an open neighborhood $W_i$ of $\BFx_i^*$ such that
\begin{equation}
    f_i(\BFx_i^*,\BFx_{-i}^*)
    \geq
    f_i(\BFx_i,\BFx_{-i}^*)
    \quad
    \text{for every }
    \BFx_i\in
    W_i\cap\widetilde{\mathcal X}_i(\BFx_{-i}^*).
    \label{eq:local_generalized_nash}
\end{equation}
\end{definition}

For brevity, we refer to a local generalized Nash equilibrium as a
local equilibrium. Definition~\ref{def:local-nash} imposes no minimum
size on the neighborhood of local optimality. Hence, a qualifying neighborhood may exclude all pricing deviations that change household
participation. Local optimality within a fixed insured portfolio
provides limited economic information because pricing incentives to
gain or lose policyholders are not evaluated. To examine economically
meaningful deviations that can alter household participation, we fix a
pricing radius $r>0$ and define
\begin{equation}
    \mathcal V_i(r;\BFx^*)
    :=
    \left\{
        \BFx_i=(\boldsymbol{\lambda}_i,A_i,M_i)
        \in\widetilde{\mathcal X}_i(\BFx_{-i}^*):
        \left\|
            \boldsymbol{\lambda}_i-\boldsymbol{\lambda}_i^*
        \right\|_2
        \leq r
    \right\}.
    \label{eq:pricing_deviation_set}
\end{equation}
For an objective tolerance $\epsilon_{\mathrm{eq}}\geq0$, a feasible
profile $\BFx^*$ satisfies the prescribed-radius approximate equilibrium
condition if
\begin{equation}
    f_i(\BFx_i^*,\BFx_{-i}^*)
    \geq
    \sup_{\BFx_i\in\mathcal V_i(r;\BFx^*)}
    f_i(\BFx_i,\BFx_{-i}^*)
    -
    \epsilon_{\mathrm{eq}},
    \qquad i=1,\ldots,N.
    \label{eq:epsilon-r-local-nash}
\end{equation}
The insolvency tolerance $\varepsilon$ determines feasibility, whereas
$\epsilon_{\mathrm{eq}}$ bounds the permitted payoff improvement.
When $\epsilon_{\mathrm{eq}}=0$, satisfaction of \eqref{eq:epsilon-r-local-nash} for any $r>0$ implies the local equilibrium condition
in Definition~\ref{def:local-nash}. The computational goal is to
identify candidate profiles and assess their equilibrium quality
against the prescribed-radius criterion within a few minutes.

\subsection{Computational Structure}
\label{subsec:computational_challenges}

Pricing and reinsurance play distinct roles in determining the insurer's capital path. With competitors' demand fixed, pricing determines household participation and the resulting insured portfolio, whereas reinsurance determines how much of the portfolio loss is retained.
The resulting premium income, claim payments, and reinsurance cash flows jointly determine financial performance and insolvency feasibility.

The effect of pricing is discontinuous because household participation
is threshold-based. As long as participation remains unchanged,
scenario-specific losses and deductibles remain fixed, while premium
income varies with prices. Crossing a purchase threshold changes the
insured portfolio discretely, producing abrupt changes in the objective
and potentially in insolvency feasibility.

\begin{figure}[htp]
\centering
\subfloat[Insolvency ratio versus $\lambda_8$]{%
    \resizebox*{7cm}{!}{%
        \includegraphics{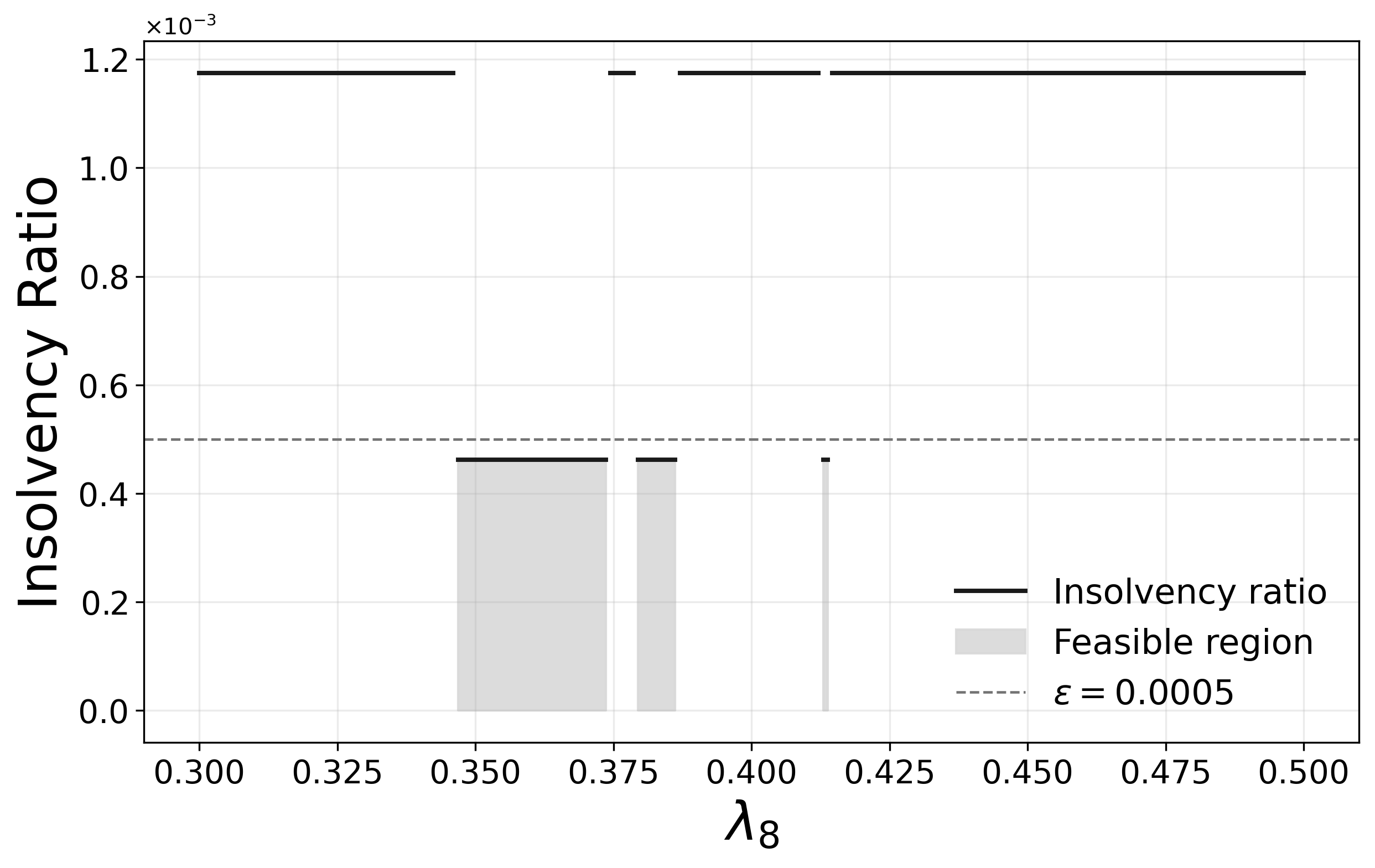}}%
    \label{fig:insolvency-price}}
\subfloat[Insolvency ratio over $(A,M)$]{%
    \resizebox*{6.3cm}{!}{%
        \includegraphics{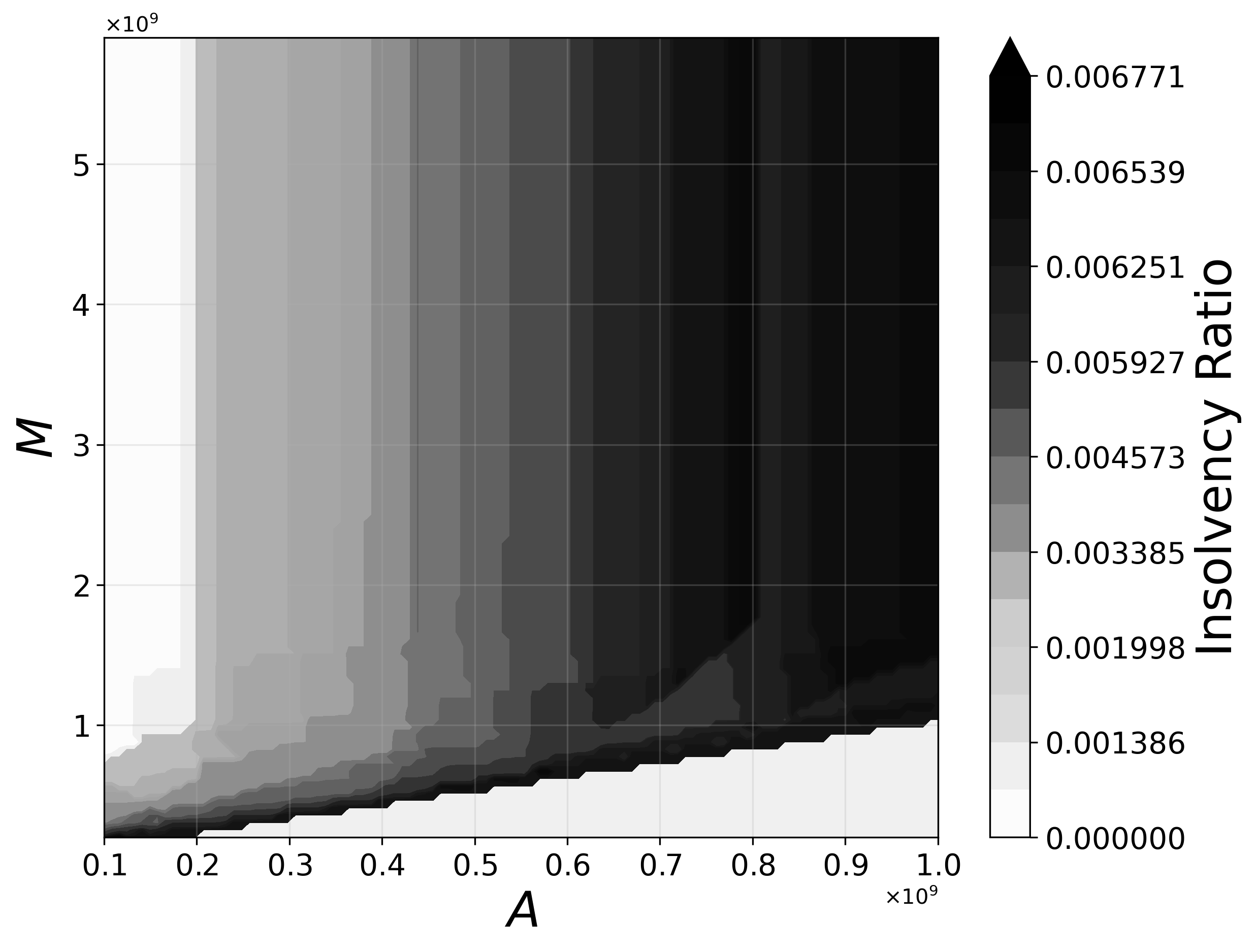}}%
    \label{fig:insolvency-reinsurance}}
\caption{
Insolvency feasibility with the remaining decisions fixed.
(a) Variation in the insolvency ratio with $\lambda_8$; the dashed line
marks the tolerance $\varepsilon=0.0005$, with values at or below the
line satisfying the insolvency constraint.
(b) Variation in the insolvency ratio across reinsurance contracts.
}
\label{fig:insolvency}
\end{figure}

Reinsurance has a different computational structure. Conditional on
the insured portfolio, excess-of-loss recoveries are continuous and
piecewise linear in $(A_i,M_i)$. Reinsurance pricing and multi-period cash flows introduce additional nonlinearities, but the original financial objective remains continuous, although potentially
nonsmooth, in the reinsurance parameters. By contrast, changes in pricing or reinsurance can affect whether the insurer is insolvent in a particular year and scenario. Because the insolvency measure aggregates
indicators of such crossings, the constraint is discontinuous and the
resulting feasible set can be nonconvex. Figure~\ref{fig:insolvency}
illustrates the feasibility structure along a pricing coordinate and
across the reinsurance parameter space.

Accordingly, to solve the response problem efficiently, we develop tailored optimization procedures for pricing and reinsurance. The next section presents the resulting computational framework in detail.

\section{Computational Framework}
\label{sec:computational_framework}

\begin{figure}[htp]
\centering
\subfloat[Objective versus $\lambda_8$]{%
    \resizebox*{7.2cm}{!}{%
        \includegraphics{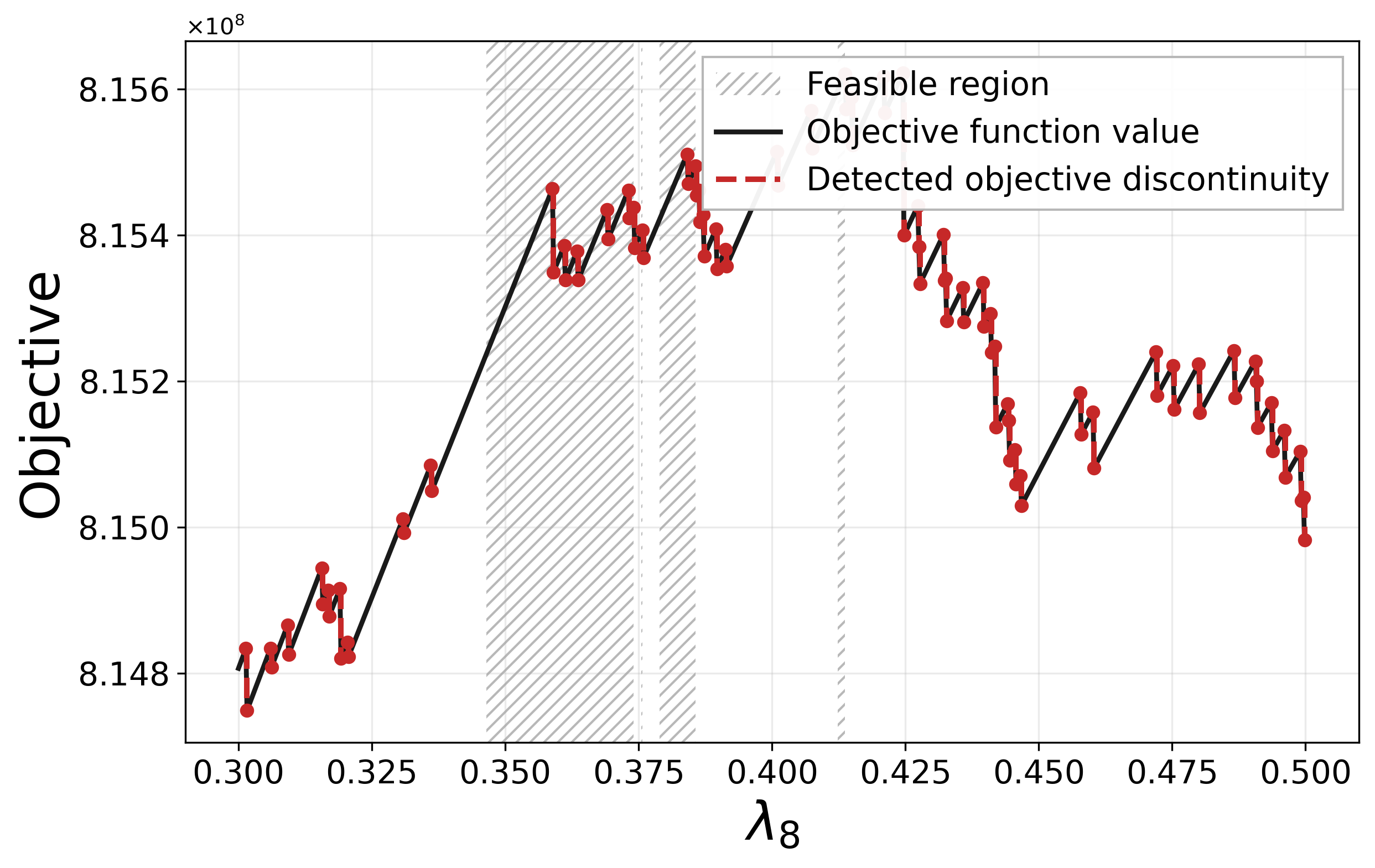}}%
    \label{fig:objective-price}}
\subfloat[Objective over $(A,M)$]{%
    \resizebox*{6cm}{!}{%
        \includegraphics{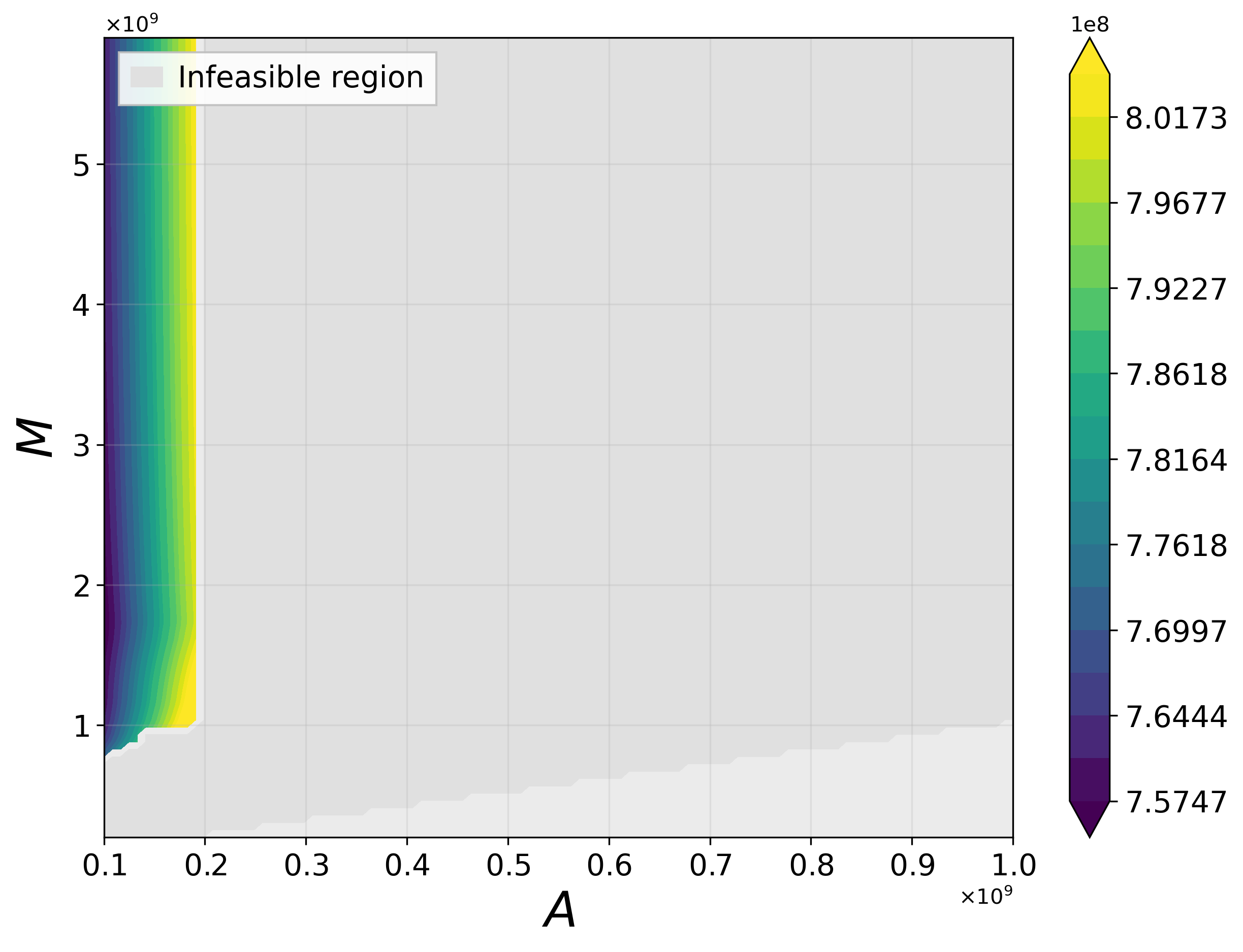}}%
    \label{fig:objective-reinsurance}}
\caption{
Objective landscapes for the insurer response problem.
(a) Variation in the objective with $\lambda_8$, including
discontinuities induced by household participation thresholds; shaded
values violate the insolvency constraint.
(b) Variation in the objective over the reinsurance parameters $(A,M)$,
with infeasible contracts masked.
}
\label{fig:surface}
\end{figure}

We develop a structured procedure for computing approximate responses to~\eqref{eq:br} with competitors' strategies fixed. The framework uses separate but coordinated procedures for pricing and reinsurance, reflecting the different objective landscapes illustrated in Figure~\ref{fig:surface}. The pricing objective exhibits abrupt changes associated with household participation thresholds, whereas the reinsurance objective is continuous but nonlinear over the admissible contracts. We first consider reinsurance optimization for fixed pricing parameters and then turn to the discontinuous pricing problem.

\subsection{Optimization of Reinsurance Parameters}

For fixed pricing parameters $\BFlambda_i$ and competitors' strategies
$\BFx_{-i}$, household participation and the insured portfolio remain
fixed. The financial objective is therefore continuous in the
reinsurance parameters $(A_i,M_i)$, although it remains nonlinear and
potentially nonsmooth. The insolvency constraint is discontinuous, and
the resulting feasible set can be nonconvex. Nevertheless, the
low-dimensional and locally structured objective landscape supports
gradient-based refinement when initialized from promising starting
solutions.

\begin{figure} [htp]
\centering
\subfloat{%
\resizebox*{6.7cm}{!}{\includegraphics{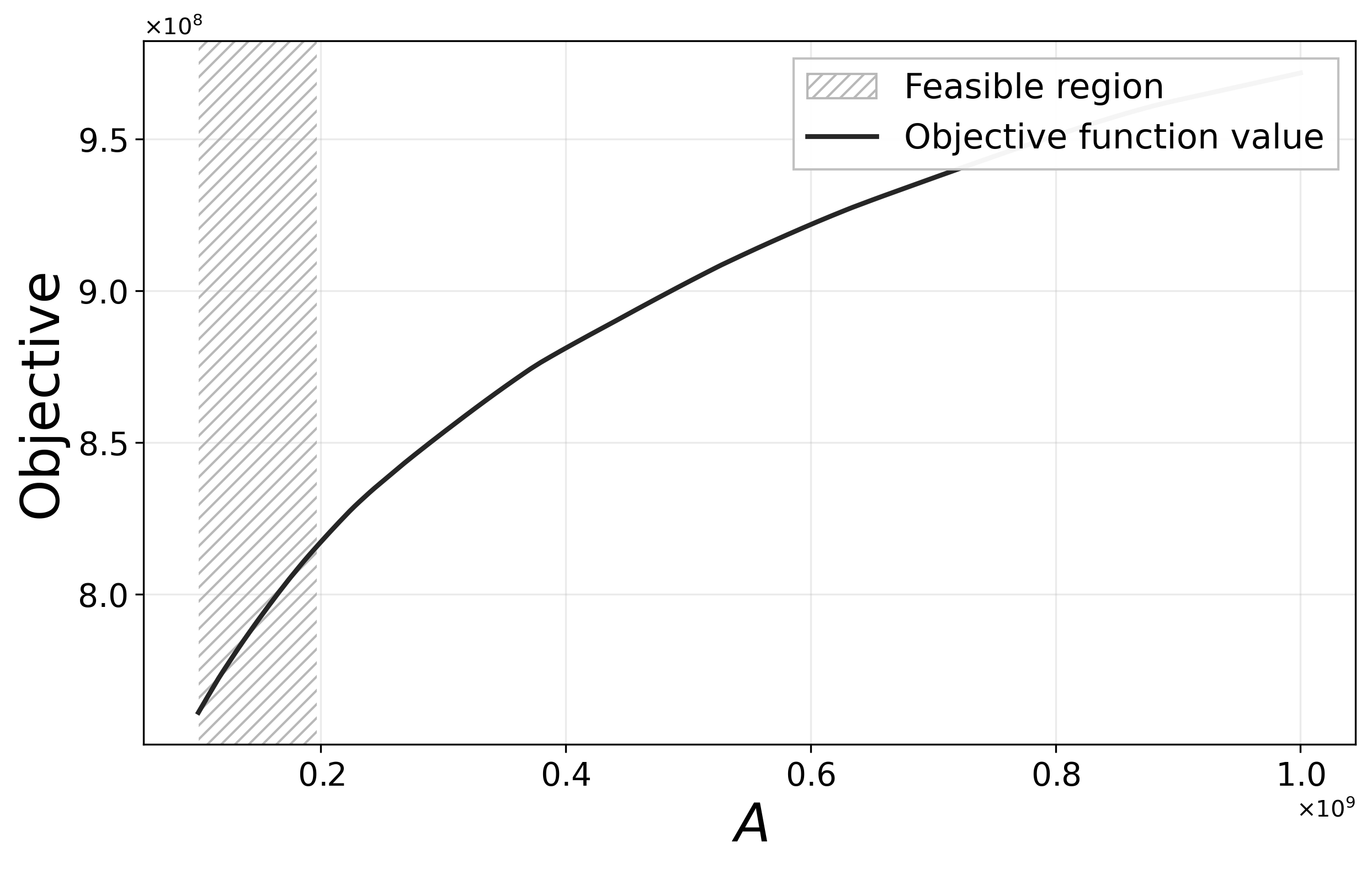}}\label{fig:A}}
\subfloat{%
\resizebox*{6.7cm}{!}{\includegraphics{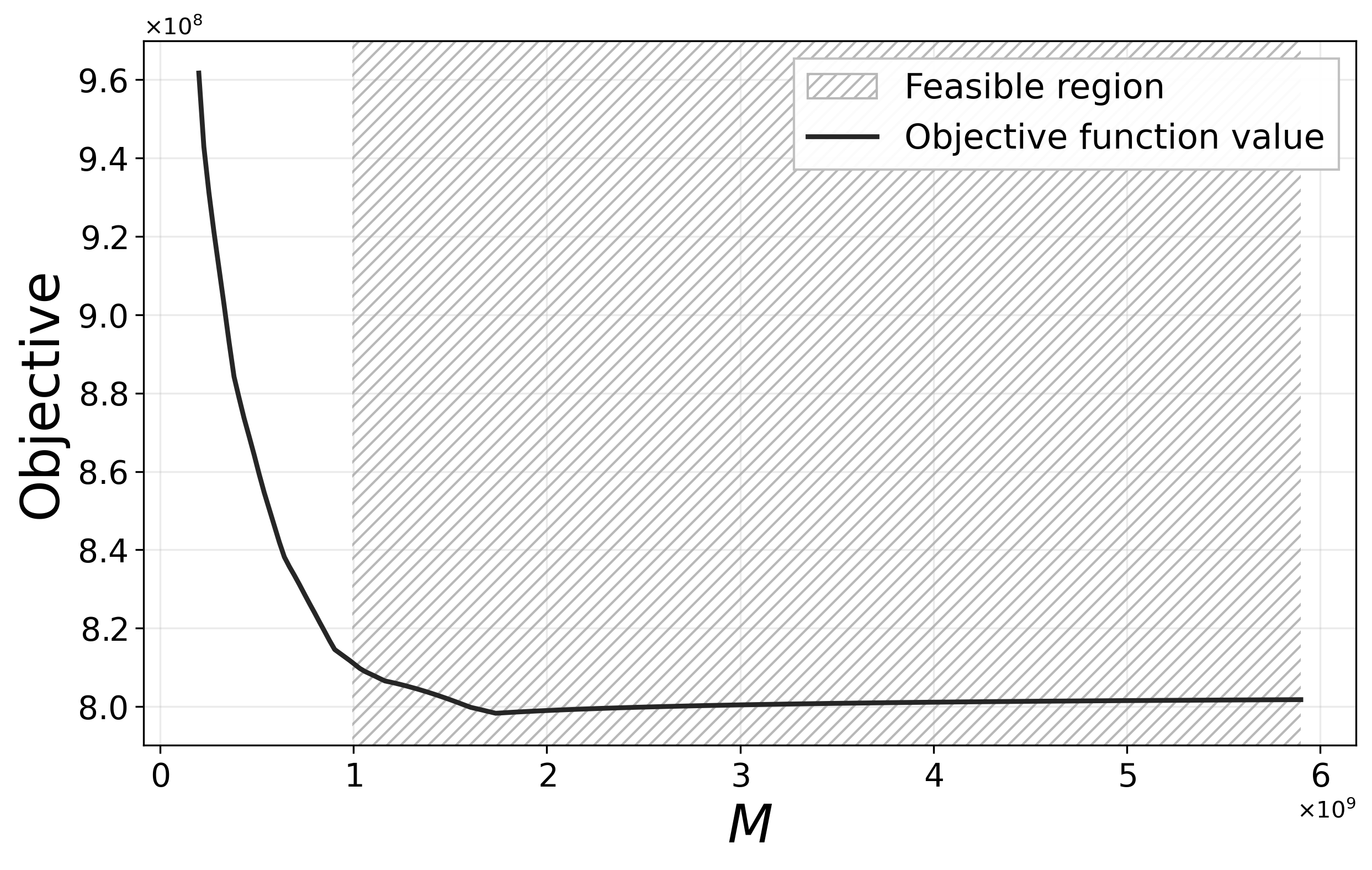}}\label{fig:M}}
\caption{
Objective function trajectory with respect to $A$ and $M$.
}
\label{fig:objective-AM}
\end{figure}

Motivated by this structure, we first perform a coarse grid search over the
reinsurance parameter space in order to identify sufficiently good candidate solutions. Starting from the best grid-search solution, we then apply a gradient-based local optimization method to further refine the reinsurance parameters. This approach provides two practical advantages. First, the coarse grid search helps to identify feasible initial solutions. Second, although no theoretical guarantee of global optimality is available, the grid-search initialization often enables the local optimization procedure to locate near-globally optimal solutions in practice. In addition, the coarse grid search can be performed efficiently using GPU parallelization, since function evaluations at different candidate solutions are computed simultaneously in batch form through vectorized Monte Carlo simulation operations. As a result, the coarse grid search does not become a computational bottleneck in practice. 

The optimization problem is solved using a penalty-augmented objective that accounts for feasibility violations. We include two penalty terms: one for violations of the insolvency constraint and another for negative demand. The latter is required because the demand attributed to competing insurers is held fixed when solving~\eqref{eq:br}. Consequently, the demand available to firm \(i\) is determined as the residual after accounting for the demand already captured by competing firms. For some candidate pricing vectors \(\BFlambda_i\), this residual may become negative in certain regions. Because negative demand has no economic meaning, we penalize such candidate solutions. Specifically, the penalty-augmented objective is defined as
\begin{equation}
\label{eq:penalty-objective}
    \Phi(\BFx_i,\BFx_{-i},\zeta_{\mathrm{ins}},\xi)
    =
    F(\BFx_i,\BFx_{-i},\xi)
    -
    \zeta_{\mathrm{ins}} \Psi_{\mathrm{ins}}(\BFx_i,\BFx_{-i},\xi)
    -
    \Psi_{\mathrm{dem}}(\BFx_i,\BFx_{-i},\xi),
\end{equation}
where \(\Psi_{\mathrm{ins}}\) denotes the insolvency constraint-related penalty term, \(\Psi_{\mathrm{dem}}\) denotes the negative-demand penalty term, and \(\zeta_{\mathrm{ins}}>0\) is a penalty parameter. Combining the two-stage search strategy with the penalty-augmented objective yields the reinsurance optimization procedure summarized in Algorithm~\ref{alg:AM-optimize}.

\begin{algorithm}[htp]
\footnotesize
\caption{Optimization of Reinsurance Parameters}
\label{alg:AM-optimize}
\begin{algorithmic}[1]
\Require Fixed pricing parameters $\BFlambda_i$, competitors' decisions $\BFx_{-i}$, step size sequence $\{\alpha_k\}$,
insolvency penalty parameter $\zeta_{\mathrm{ins}}$, the maximum number of iterations $K$, 
and feasible reinsurance region $\mathcal{R}$.

\State Perform a coarse grid search over the reinsurance parameter space and set
\[
    (A_0,M_0)
    =
    \argmax_{(A,M)\in\mcG}
    \sum_{j\in\mcS}
    p_j
    \Phi(\BFlambda_i,A,M,\BFx_{-i},\zeta_{\mathrm{ins}},\xi_j),
\]
where $\mcG$ denotes the grid-search candidate set.

\State Set $k=0$ and $\alpha_k=\alpha_0$.

\For{$k=0,1,2,\ldots, K-1$}

\State Evaluate
\(
    \bar{\Phi}_k
    =
    \sum_{j\in\mcS}
    p_j
    \Phi(\BFlambda_i,A_k,M_k,\BFx_{-i},\zeta_{\mathrm{ins}},\xi_j)
\)
and 
\(
    \BFg_k
    =
    \nabla_{(A,M)}
    \bar{\Phi}_k .
\)

\State Compute the projected candidate
\(
    (A_k^{s},M_k^{s})
    =
    \Pi_{\mathcal R}\!\left(
        (A_k,M_k)+\alpha_k\BFg_k
    \right).
\)
\State Evaluate
\(
    \bar{\Phi}_k^{s}
    =
    \sum_{j\in\mcS}
    p_j
    \Phi(\BFlambda_i,A_k^{s},M_k^{s},\BFx_{-i},\zeta_{\mathrm{ins}},\xi_j).
\)

\If{$\bar{\Phi}_k^{s} \ge \bar{\Phi}_k$}
    \State Accept the update and set
    \(
        (A_{k+1},M_{k+1})
        =
        (A_k^{s},M_k^{s}).
    \)
\Else
    \State Reject the update and set
    \(
        (A_{k+1},M_{k+1})
        =
        (A_k,M_k).
    \)
\EndIf
\EndFor
\State \Return [$A_K, M_K$]
\end{algorithmic}
\end{algorithm}

For fixed pricing decisions \(\BFlambda_i\) and competitors' decisions
\(\BFx_{-i}\), Algorithm~\ref{alg:AM-optimize} approximates the
reinsurance optimization operator
\(
    \mathcal P(\BFlambda_i;\BFx_{-i})
\),
defined by
\[
    \mathcal P(\BFlambda_i;\BFx_{-i})
    \in
    \arg\max_{(A,M)\in\mathcal R}
    f(\BFlambda_i,A,M,\BFx_{-i})
    \quad
    \text{s.t.}
    \quad
    g(\BFlambda_i,A,M,\BFx_{-i})
    \leq\varepsilon.
\]
The original problem~\eqref{eq:br} can then be reduced to
\begin{equation}
    \max_{\BFlambda_i\in
        \Lambda_i^{\mathrm{feas}}(\BFx_{-i})}
    f\!\left(
        \BFlambda_i,
        \mathcal P(\BFlambda_i;\BFx_{-i}),
        \BFx_{-i}
    \right),
    \label{eq:projected-br}
\end{equation}
where \(\Lambda_i^{\mathrm{feas}}(\BFx_{-i})\) is the set of pricing
vectors in \(\Lambda\) that admit at least one reinsurance contract in
\(\mathcal R\) satisfying the insolvency constraint.

\begin{figure} [htp]
\centering
\subfloat[Objective trajectory vs. $\lambda_1$]{%
\resizebox*{6.7cm}{!}{\includegraphics{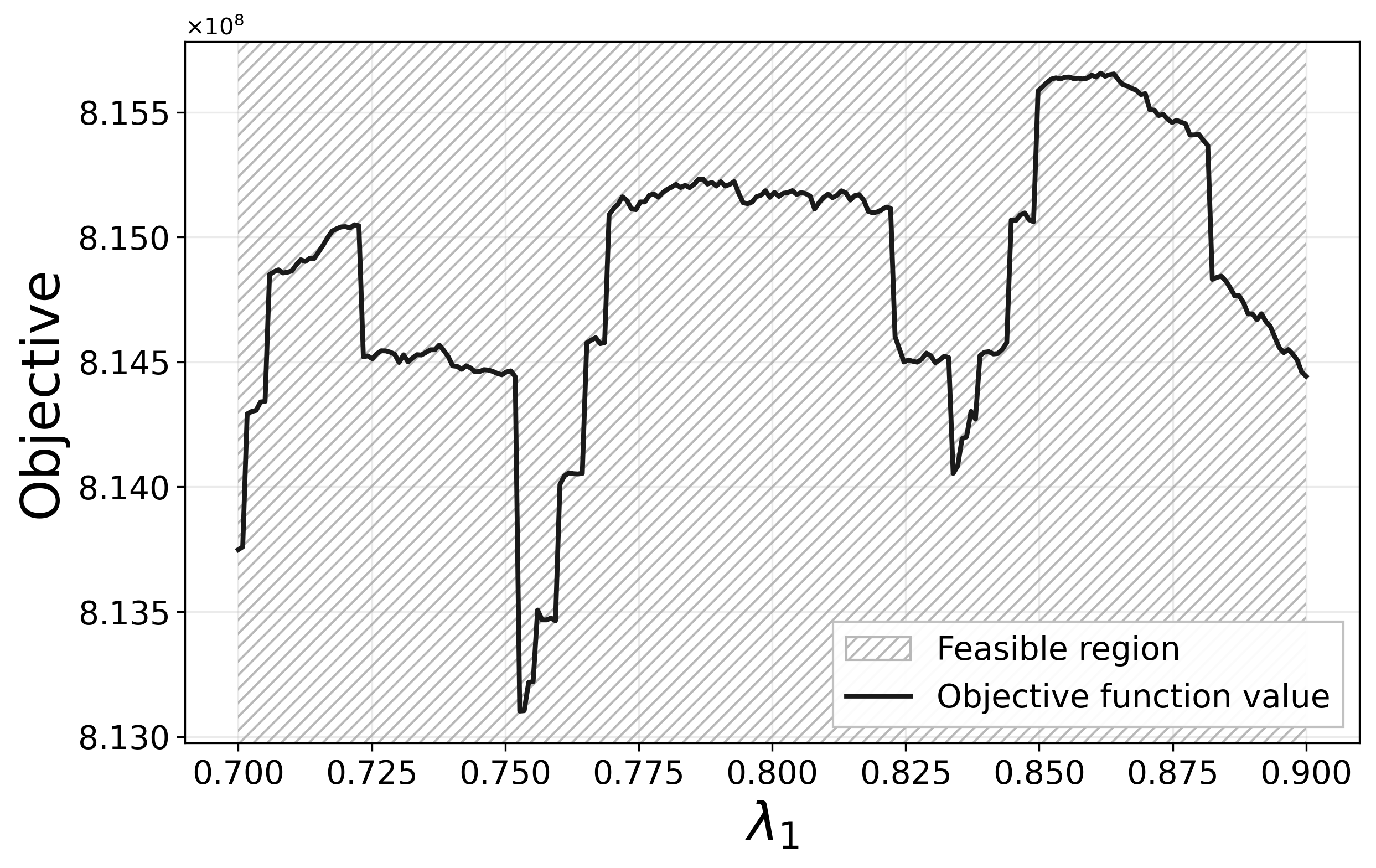}}\label{fig:4-1}}
\subfloat[Objective trajectory vs. $\lambda_8$]{%
\resizebox*{6.7cm}{!}{\includegraphics{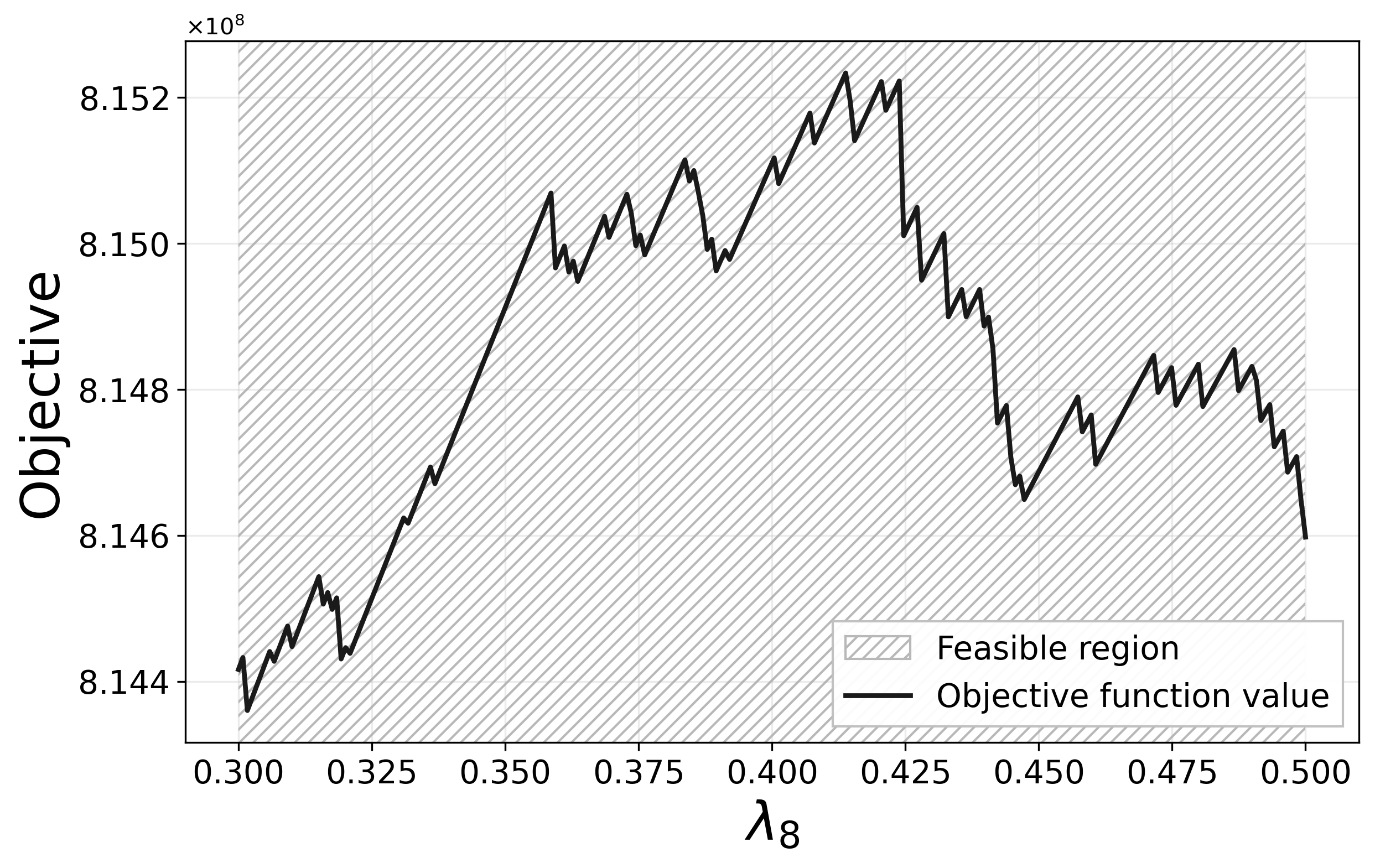}}\label{fig:4-2}}
\caption{
Objective function trajectories of~\eqref{eq:projected-br} with respect to $\lambda_1$ and $\lambda_8$, highlighting the heterogeneous structure of the projected pricing optimization problem across regions.
}
\label{fig:surface-reopt}
\end{figure}

As a result, the pricing optimization procedure can be viewed as searching over pricing decisions while optimizing the associated reinsurance parameters through the operator \(\mathcal P\). Figure~\ref{fig:surface-reopt} further illustrates that, after applying \(\mathcal P\), feasible reinsurance decisions can be found for a wide range of pricing vectors. Thus, the problem~\eqref{eq:projected-br} effectively removes much of the complexity associated with the insolvency constraint, at the expense of introducing additional nonsmoothness into the objective function.

\subsection{Optimization of Pricing Parameters}

\begin{figure} [htp]
\centering
\subfloat[$\mu = 0.001$]{%
\resizebox*{6.7cm}{!}{\includegraphics{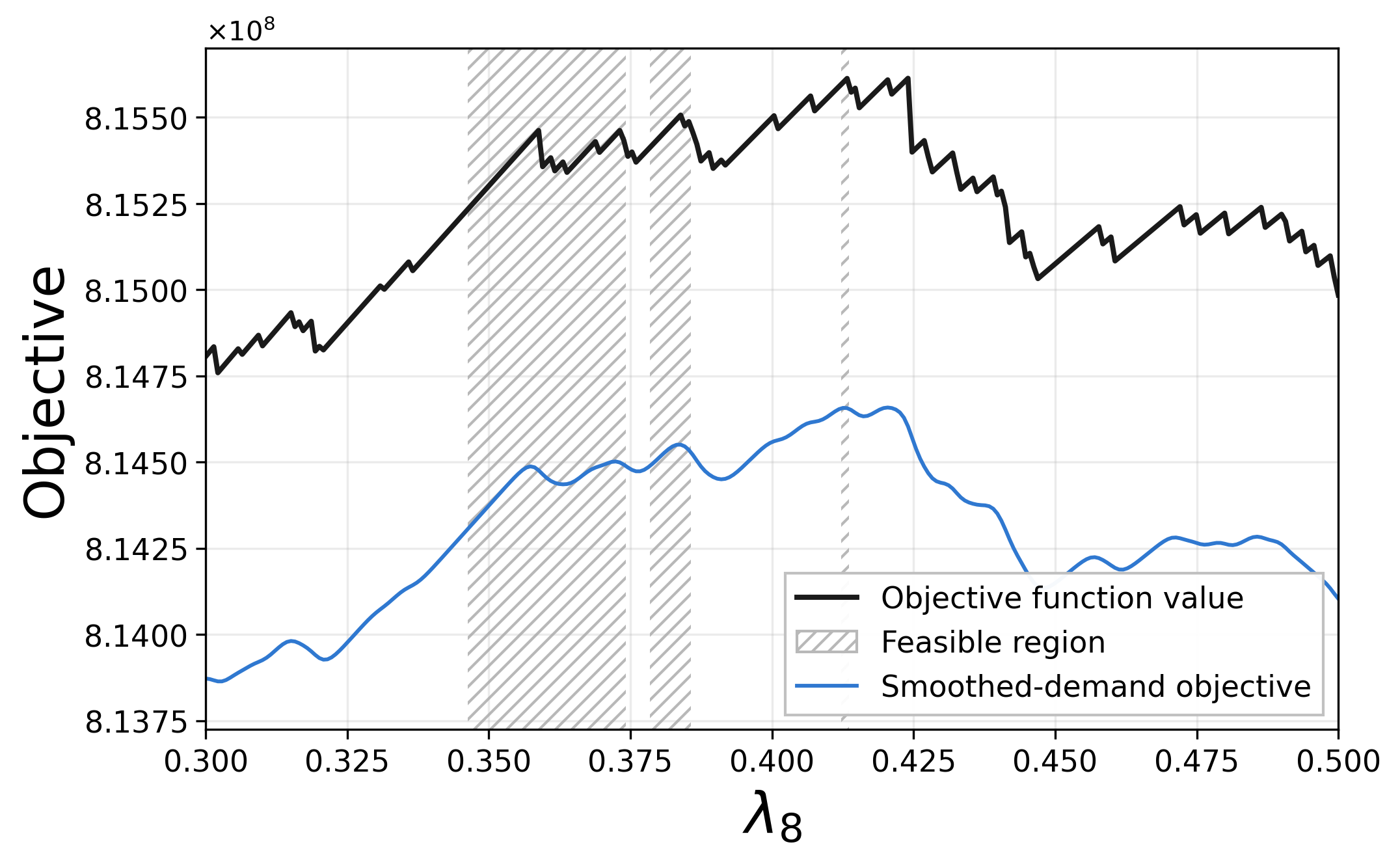}}\label{fig:mu=0.001}}
\subfloat[$\mu = 0.0001$]{%
\resizebox*{6.7cm}{!}{\includegraphics{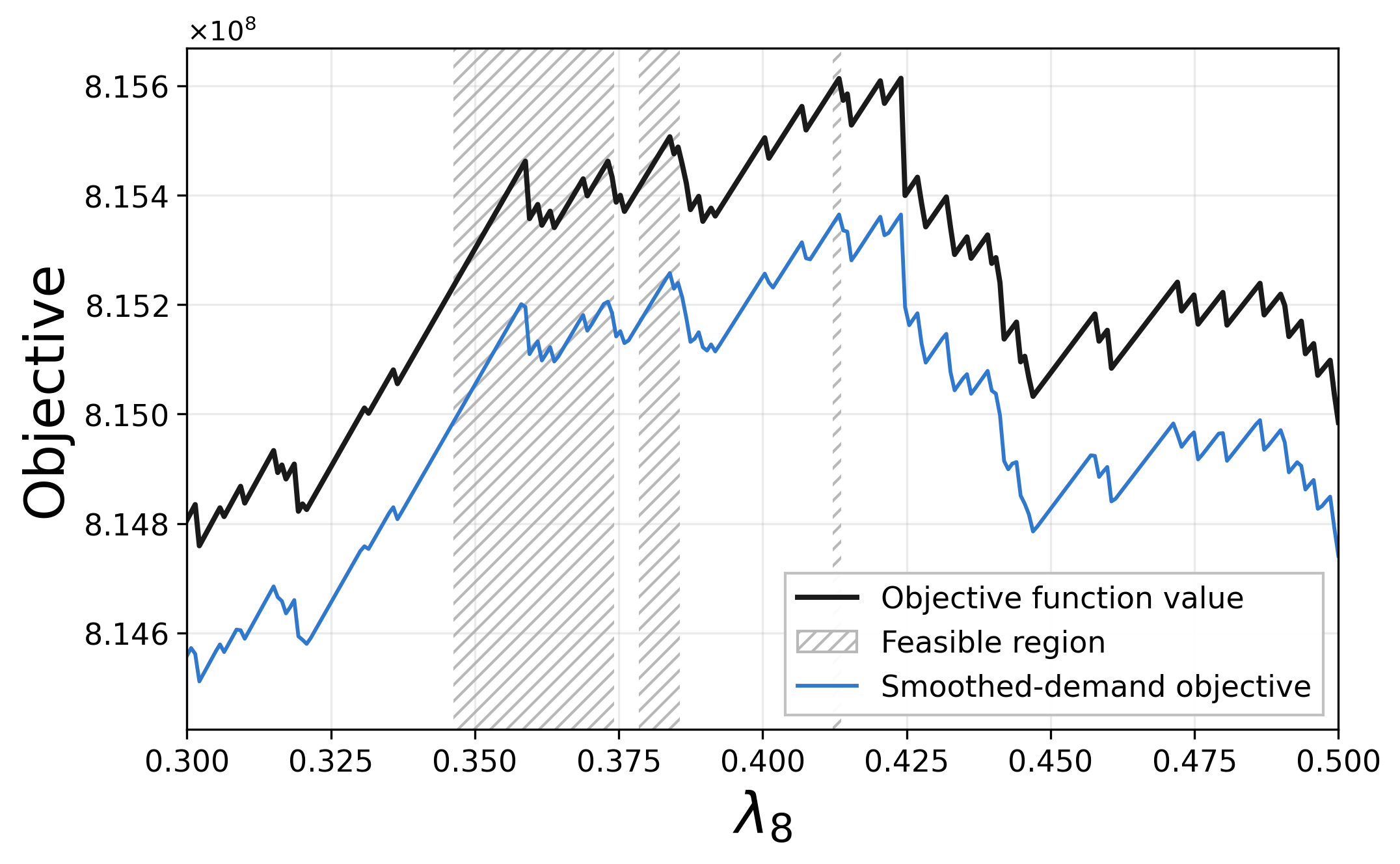}}\label{fig:mu=0.0001}}
\caption{
Objective function trajectory with respect to $\lambda_8$ under different smoothing parameters.
}
\label{fig:smoothed-demand}
\end{figure}

In this section, we focus on the pricing parameters. In particular, given fixed reinsurance parameters $(A,M)$ and competitors’ decisions $\BFx_{-i}$, the representative insurer solves problem~\eqref{eq:br} with respect to its pricing vector $\BFlambda$. 
To address discontinuities in the objective function induced by threshold-based demand responses, we introduce a smoothing approximation of the demand function. Specifically, we approximate the demand indicator using a sigmoid-based formulation:
\(
\mathbf{1}\{\lambda \in [\underline{\lambda}, \overline{\lambda}]\} 
\approx 
S_{\mu}(\lambda - \underline{\lambda}) - S_{\mu}(\lambda - \overline{\lambda}),
\)
where 
\(
S_{\mu}(x) = (1 + \exp(-x/\mu))^{-1}
\)
and $\mu > 0$ controls the smoothness of the approximation; as $\mu \to 0$, the function converges to the original demand function. Here, $\underline{\lambda}$ and $\overline{\lambda}$ denote the lower and upper threshold values that determine homeowners' insurance purchase decisions. Let \(F_\mu(\BFx_i,\BFx_{-i},\xi)\) and \(\Phi_\mu(\BFx_i,\BFx_{-i},\zeta_{\mathrm{ins}},\xi)\) denote the smoothed versions of \(F\) and \(\Phi\), respectively, obtained by replacing the demand indicator with its sigmoid approximation.
Although min and max operations leave nonsmooth kinks in the formulation, the smoothed objective is differentiable almost everywhere, enabling efficient derivative computation with standard automatic differentiation tools in practice~\citep{paszke2019pytorch}.
However, a fundamental trade-off arises in practice. Larger values of $\mu$ improve numerical stability but may significantly distort the structure of the underlying problem. In contrast, smaller values of $\mu$, while providing a closer approximation to the original discontinuous function, can lead to numerically unstable or ill-conditioned gradients, resulting in unreliable updates; see Figure~\ref{fig:smoothed-demand}.

To address these challenges, we adopt a hybrid optimization approach that combines gradient-based methods with derivative-free methods. Optimization is carried out within a trust-region framework, where a local quadratic model is constructed using gradient information from \(F_{\mu}(\BFx_i, \BFx_{-i}, \xi)\) when numerically reliable, and finite-difference (FD) approximations based on evaluations of the original objective otherwise; see \ref{appdx:s3} for details on the construction of the FD-based gradient and Hessian approximations.

\begin{algorithm}[htp]
\footnotesize
\caption{Optimization of Pricing Parameters}
\label{alg:lambda-optimize}
\begin{algorithmic}[1]
\Require Current pricing vector $\BFlambda_i^k$,
reinsurance parameters $(A,M)$, competitors' decisions $\BFx_{-i}^k$,
radii $\Delta_k,\Delta_{\max}$,
thresholds $\eta\in(0,1)$ and $\Delta_{\mathrm{th}}>0$,
smoothing parameter $\mu>0$,
update factors $\gamma_1>1>\gamma_2>0$,
and improvement tolerance $\tau_{\mathrm{acc}}>0$.

\State Define
\[
\bar{\Phi}^{(k)}(\BFlambda;a,m)
:=
\begin{cases}
    \sum_{j\in\mcS}p_j\Phi_\mu(\BFlambda,\zeta_{\mathrm{ins}}^k,\xi_j;
             a,m,\BFlambda_{-i}),
    & \Delta_k>\Delta_{\mathrm{th}},\\
    \sum_{j\in\mcS}p_j\Phi(\BFlambda,\zeta_{\mathrm{ins}}^k,\xi_j;
         a,m,\BFlambda_{-i}),
    & \text{otherwise}.
\end{cases}
\]

\State Evaluate
$\bar{\Phi}_k=\bar{\Phi}^{(k)}(\BFlambda_i^k;A,M)$.

\If{$\Delta_k>\Delta_{\mathrm{th}}$}
    \State Using automatic differentiation, compute
    \(
        \BFg_k
        =
        \nabla_{\BFlambda}\bar{\Phi}_k,
    \)
    and
    \(
        \sfH_k
        =
        \nabla_{\BFlambda}^2\bar{\Phi}_k.
    \)
    
    \State Set the direct-search candidate set $\mcC_k=\emptyset$.
\Else
    \State Compute FD approximations $\BFg_k,\sfH_k$ using the stencil points.
    \State Form $\mcC_k$ from the stencil and additional exploration points, totaling 33 points.
\EndIf

\State Construct the quadratic model 
$M_k(\BFs)=\bar{\Phi}_k+\BFs^\intercal\BFg_k
+\frac12\BFs^\intercal\sfH_k\BFs$.

\State Compute the trust-region candidate 
$\BFs^k
\in
\argmax_{\substack{
    \|\BFs\|\le\Delta_k\\
    \BFlambda_i^k+\BFs\in\Lambda
}}
M_k(\BFs)$
and set
$\widetilde{\BFlambda}_i^k=\BFlambda_i^k+\BFs^k$.

\State Evaluate the candidates using $\bar{\Phi}^{(k)}(\cdot;A,M)$, and select the best candidate
\[
\widehat{\BFlambda}_i^k
\in
\argmax_{\BFlambda\in
\{\widetilde{\BFlambda}_i^k\}\cup\mcC_k}
\bar{\Phi}^{(k)}(\BFlambda;A,M).
\]

\State Given $\widehat{\BFlambda}_i^k$ and $\BFx_{-i}^k$,
obtain $(\widehat A,\widehat M)$ using
Algorithm~\ref{alg:AM-optimize} with objective
$\bar{\Phi}^{(k)}(\widehat{\BFlambda}_i^k;\cdot,\cdot)$.

\State Evaluate
$\bar{\Phi}_k^{\mathrm{cand}}
=\bar{\Phi}^{(k)}
(\widehat{\BFlambda}_i^k;\widehat A,\widehat M)$
and compute
\[
\rho_k=
\frac{\bar{\Phi}_k^{\mathrm{cand}}-\bar{\Phi}_k}
{M_k(\widehat{\BFlambda}_i^k-\BFlambda_i^k)-M_k(\mathbf0)}.
\]

\If{$\rho_k\ge\eta$ and
$\bar{\Phi}_k^{\mathrm{cand}}-\bar{\Phi}_k
>\tau_{\mathrm{acc}}\max\{1,|\bar{\Phi}_k|\}$}
    \State Set
    $\BFlambda_i^{k+1}=\widehat{\BFlambda}_i^k$ and
    $\Delta_{k+1}=\min\{\gamma_1\Delta_k,\Delta_{\max}\}$.
\Else
    \State Set
    $\BFlambda_i^{k+1}=\BFlambda_i^k$ and
    $\Delta_{k+1}=\gamma_2\Delta_k$.
\EndIf

\State \Return $[\BFlambda_i^{k+1},\Delta_{k+1}]$.
\end{algorithmic}
\end{algorithm}

When FD approximations are required, the resulting function evaluations can also be utilized for direct-search exploration. In particular, constructing FD-based gradient and Hessian approximations already requires evaluating the objective at multiple nearby points within the trust-region neighborhood. These stencil points can therefore be reused as candidate direct-search steps without additional function evaluations. To broaden the search beyond the coordinate directions, we also evaluate additional points generated in both directions along each vector of a randomly generated orthonormal basis, with the points projected onto the pricing box as needed. The combined set contains 33 points, including the current iterate. The best candidate among these points and the trust-region candidate is then evaluated after reoptimizing the reinsurance parameters and accepted if the acceptance criteria are satisfied. The resulting optimization procedure is summarized in Algorithm~\ref{alg:lambda-optimize}.

\begin{remark}[Role of \(\tau_{\mathrm{acc}}\) in Algorithm~\ref{alg:lambda-optimize}]
Direct-search candidates are not constructed to satisfy the sufficient model decrease condition associated with the Cauchy step in standard trust-region methods. We therefore require a minimum relative objective improvement determined by \(\tau_{\mathrm{acc}}\), preventing the acceptance of arbitrarily small improvements. In the \((\epsilon,r)\) local Nash condition of Section~\ref{sec:intro-LNE}, \(\tau_{\mathrm{acc}}\) controls the practical tolerance associated with \(\epsilon\). 
\end{remark}

\subsection{Finding Equilibria}

To compute equilibrium solutions, we exploit the homogeneity of firms and restrict attention to symmetric equilibrium solutions satisfying
\(
\BFx_i=\BFx_j
\)
for all \(i\) and \(j\). Consequently, it suffices to optimize the strategy of a representative firm while treating the demand induced by competing firms as fixed. In particular, we repeatedly compute a better-response strategy for the representative firm using Algorithms~\ref{alg:AM-optimize} and~\ref{alg:lambda-optimize}. The resulting strategy is then assumed to be adopted by all firms, and the corresponding competitors' demand estimates are updated for the next iteration. This process is repeated until the strategy profile stabilizes. 
However, directly updating the competitors' demand estimates using the newly computed solution \((\BFlambda_i^{k+1})\) can lead to numerical instability. In particular, large changes in the pricing parameters may induce abrupt changes in the effective demand, causing negative-demand regions to suddenly appear or disappear during the equilibrium iterations. As a result, both the competitors' demand estimates and the resulting pricing updates may fluctuate significantly across successive iterations, potentially leading to cycling behavior.

To improve numerical stability, we therefore apply a damped update rule to the competitors' demand estimates during the outer equilibrium iterations. Rather than fully replacing the previous demand estimates using the newly computed solution, we update the competitors' demand estimates gradually according to
\(
    D_{k+1}^{\mathrm{other}}
    =
    (1-\alpha)D_k^{\mathrm{other}}
    +
    \alpha \widetilde{D}_k^{\mathrm{other}},
\)
where \(D_k^{\mathrm{other}}\) denotes the current estimate of the demand captured by competing firms, \(\widetilde{D}_k^{\mathrm{other}}\) denotes the updated demand estimate induced by the newly computed strategy \((\BFlambda_i^{k+1})\), and \(\alpha\in(0,1]\) is a relaxation parameter. In our numerical experiments, we found that this damped update procedure plays a crucial role in improving numerical stability and preventing cycling behavior during the equilibrium iterations. As the pricing strategies converge, the corresponding competitors' demand estimates also stabilize under the damped update procedure. Finally, we summarize the overall equilibrium computation procedure in Algorithm~\ref{alg:equilibria}.

\begin{algorithm}[htp]
\footnotesize
\caption{Finding Equilibria}
\label{alg:equilibria}
\begin{algorithmic}[1]
\Require Initial strategy profile $\BFx_i^{0}$, relaxation parameter $\alpha\in(0,1]$,
tolerance $\epsilon_{\mathrm{eq}}>0$, initial trust-region radius $\Delta_0$, maximum number of outer iterations $K_{\max}$, and the number of competing firms $N_{\mathrm{other}}.$

\State Initialize the aggregate demand estimate of competing firms,
$D_0^{\mathrm{other}}$.

\For{$k=0,1,2,\ldots,K_{\max}$}

    \State \parbox[t]{0.95\linewidth}{
    Given competitor's strategy $D_k^{\mathrm{other}}$, $A_i^{k},M_i^k$, and trust-region radius $\Delta_k$, solve the pricing subproblem using
    Algorithm~\ref{alg:lambda-optimize} and obtain
    $\BFlambda_i^{k+1}$ and $\Delta_{k+1}$.
    }

    \State \parbox[t]{0.88\linewidth}{
    Given $\BFlambda_i^{k+1}$ and $[D_k^{\mathrm{other}},A_i^{k},M_i^k]$, solve the
    reinsurance subproblem using
    Algorithm~\ref{alg:AM-optimize} and obtain
    $(A_i^{k+1},M_i^{k+1})$.
    }


    \State Compute the induced demand vector $d_k$ from
    $\BFlambda_i^{k+1}$ and set
    \(
        \widetilde{D}_k^{\mathrm{other}}
        =
        N_{\mathrm{other}}d_k.
    \)

    \State Update the aggregate demand estimate of competing firms:
    \[
        D_{k+1}^{\mathrm{other}}
        =
        (1-\alpha)D_k^{\mathrm{other}}
        +
        \alpha \widetilde{D}_k^{\mathrm{other}}.
        \label{algstep:damped-update}
    \]

    \If{$|D_{k+1}^{\mathrm{other}}-D_{k}^{\mathrm{other}}|
    \le
    \epsilon_{\mathrm{eq}}$ and $\Delta_{k+1} < \epsilon_{\mathrm{eq}}$}
        \State \textbf{Break}
    \EndIf
\EndFor
\State \Return$\BFx_i^{(k+1)}$.
\end{algorithmic}
\end{algorithm}

\section{Numerical Experiments}
\label{sec:num-exp}

We consider a competitive market consisting of four homogeneous insurers operating across eight geographical regions in North Carolina; see Figure~3 in~\citet{liu2025computing}.
Each insurer determines a vector of region specific pricing parameters
\(
\BFlambda\in\mathbb{R}^8
\)
together with two reinsurance parameters \(A\) and \(M\).
Equilibrium iterations are initialized from randomly generated pricing vectors together with their corresponding optimized reinsurance parameters \(A\) and \(M\).
Since the optimization landscape is nonconvex and discontinuous, different initializations may lead to different equilibrium solutions.
To investigate the multiplicity of equilibria, we repeat the equilibrium computation procedure from multiple random initializations.

\medskip
\noindent
\textbf{Experimental setting}
The homeowner property values, risk aversion coefficients, and hurricane loss estimates are obtained from the data set developed by \citet{liu2025computing}.
The data contain 97 hurricane event types and event specific probabilities.
Insurer performance is evaluated using 97 probability weighted scenario paths, each representing a 30 year sequence of hurricane events.
Each scenario may contain up to five hurricane events in a given year.
The annual discount factor is set to \(\beta=0.95\).

Homeowner demand is determined using a risk aversion utility model with homeowner-specific risk aversion coefficients. For homeowner \(j\) in region \(z\) and coverage type \(r\in\{\mathrm{wind},\mathrm{flood}\}\), let \(L_{j,h}^r\) denote the
loss under hurricane event \(h\), and let
\(
    B_{j,h}^r
    =
    \min\{L_{j,h}^r,2500\}
\)
denote the portion of the loss retained by the homeowner under the
deductible. Then the expected insured loss becomes
\(
    \ell_j^r
    =
    \sum_{h=1}^{97}
    p_h
    \left(
        L_{j,h}^r-B_{j,h}^r
    \right),
\)
where the expectation is taken over the 97 hurricane event types and \(p_h\) denotes the annual occurrence probability of event type \(h\). Accordingly, the annual premium charged to homeowner \(j\) in region \(z\) for coverage type \(r\) is
\(
    \pi_j^r(\lambda_z)
    =
    (1.35+\lambda_z)\ell_j^r,
\)
where \(1.35\) is the baseline premium multiplier and \(\lambda_z\) is
the insurer's pricing decision for region \(z\). A homeowner purchases coverage when the premium does not exceed the maximum willingness to pay implied by the utility model. In the experiments with an affordability constraint, the premium is additionally restricted to be no greater than \(5\%\) of the homeowner's property value. The detailed willingness-to-pay calculation and demand construction are provided in \ref{appdx:s1}.

For each candidate pricing decision, the insurer's initial capital is set equal to three times the total annual premium income collected from its insured portfolio across all eight regions.
The insolvency measure is computed by averaging the probability of negative capital across the 30 years, and its tolerance $\varepsilon$ is set to
\(5\times10^{-4}\).
Additional details concerning the reinsurance pricing and capital dynamics are given in \ref{appdx:s2}. 

\medskip
\noindent
\textbf{Computational setting}
The regional pricing parameters are restricted to
\(\lambda_z\in[0.01,5]\), while the reinsurance parameters satisfy
\(A\in[10^8,10^9]\) and
\(M\in[2\times10^8,5.9\times10^9]\).
Unless otherwise stated, the smoothing parameter is set to
\(\mu=0.1\) and $\tau_{\text{acc}} = 10^{-7}$. \ref{appdx:s5} examines the sensitivity of the computed equilibrium solutions to the choice of \(\mu\).
For each complete reoptimization of the reinsurance parameters, we first evaluate a \(100\times200\) uniform grid over the feasible ranges of \(A\) and \(M\), and then perform five gradient-ascent iterations starting from the best grid point.
For the pricing optimization, the initial trust region radius is set to \(0.5\), with a maximum radius of \(5.0\) and an acceptance threshold $\eta$ of \(0.1\).
The aggregate demand assigned to the other insurers is updated using a damping parameter of \(0.3\).

All experiments are implemented in Python using PyTorch and conducted on a computer equipped with a single NVIDIA RTX A1000 GPU with 8 GB of memory. 
To efficiently handle the large number of Monte Carlo evaluations required during the equilibrium computations, all scenario simulations, grid-search evaluations, and batched objective computations are fully vectorized and executed in parallel on the GPU. 
In particular, candidate pricing and reinsurance solutions are evaluated simultaneously across scenarios and time periods using tensor-based operations, which substantially reduces the computational overhead associated with repeated simulation-based optimization.

\subsection{Demand Functions}

\begin{figure} [htp]
\centering
\subfloat[Without an affordability cap]{%
\resizebox*{6.7cm}{!}{\includegraphics{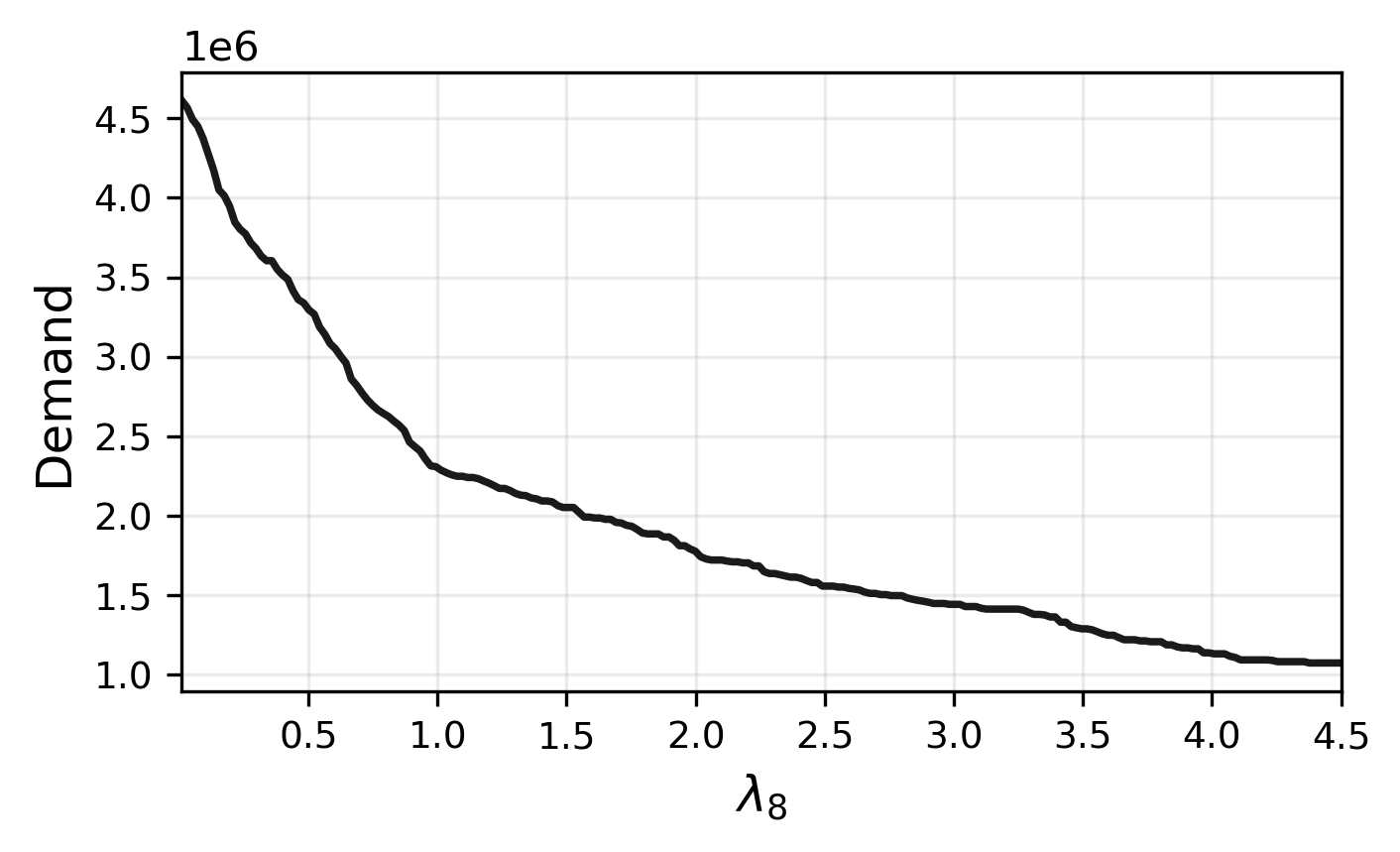}}\label{fig:demand-alphainf}}
\subfloat[With an affordability cap]{%
\resizebox*{6.7cm}{!}{\includegraphics{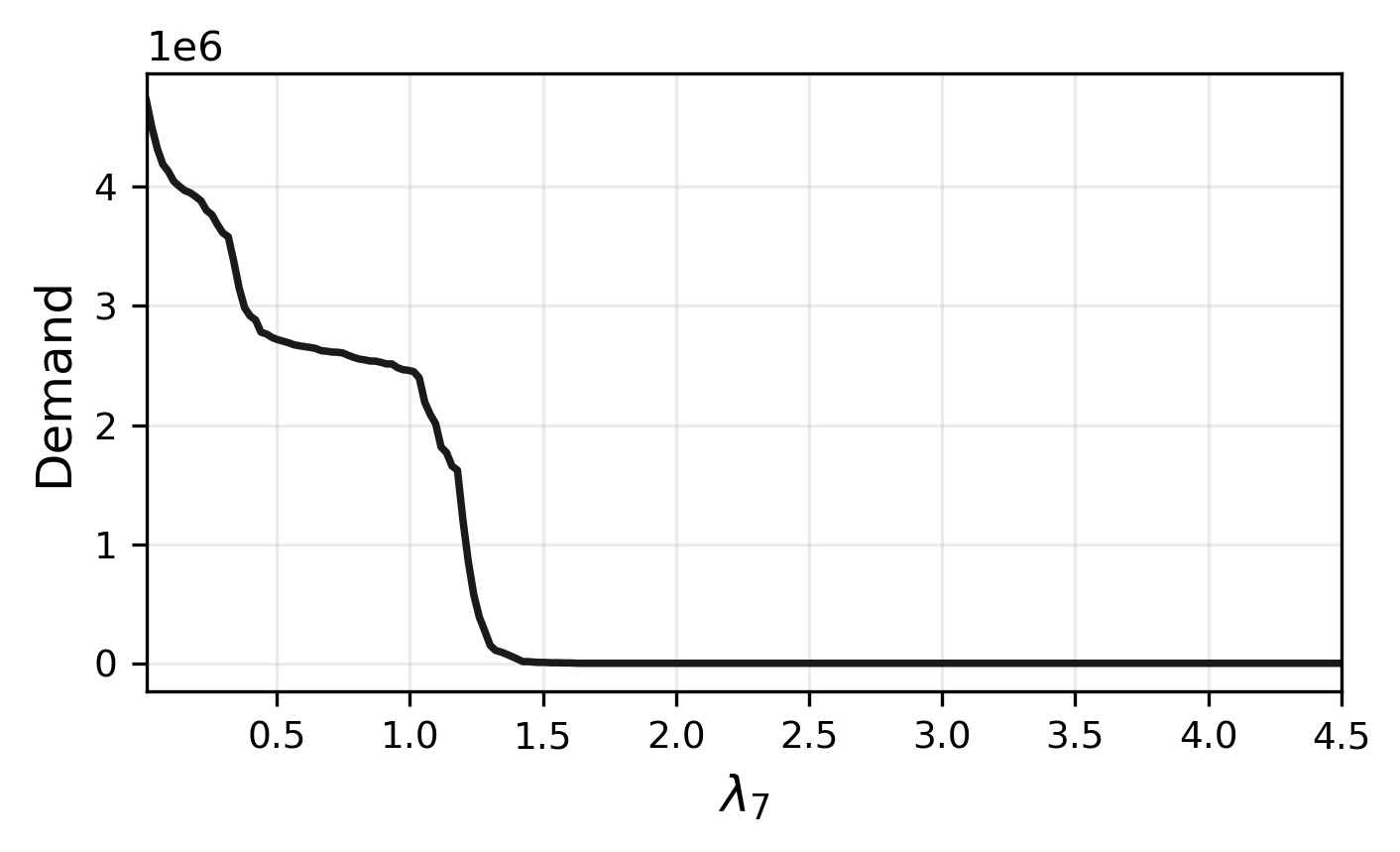}}\label{fig:demand-alpha0.05}}
\caption{
Demand as a function of the pricing parameter $\lambda_8$ under two affordability specifications:
(a) without an affordability cap ($\alpha=\infty$ in~\eqref{eq:affordability}) and
(b) with an affordability cap ($\alpha=0.05$ in~\eqref{eq:affordability}).
}
\label{fig:demand-functions}
\end{figure}

To better understand the problem structure, we begin by illustrating the aggregate demand functions induced by homeowners' threshold-based purchase decisions. Without an affordability cap, aggregate demand decreases gradually as the pricing parameter increases; see Figure~\ref{fig:demand-alphainf}. With the affordability cap, demand declines sharply over a relatively narrow range of prices and becomes nearly zero thereafter; see Figure~\ref{fig:demand-alpha0.05}. In both cases, aggregate demand is nonincreasing and nonlinear because homeowners leave the market as premiums exceed their purchase thresholds.

\subsection{Multiple Equilibrium Solutions}
Optimization problems involving insolvency-risk constraints often benefit from carefully chosen initial solutions, as appropriate initialization can substantially improve the efficiency and reliability of the optimization procedure~\citep{meesena2025safe}. Motivated by this observation, we first present an approximate local equilibrium obtained from a specific class of initial pricing vectors. We then investigate the existence of additional equilibrium solutions by initializing the algorithm from a broader collection of pricing vectors. 

\subsubsection{Without Affordability Cap}

We first consider the baseline model without an affordability cap. Table~\ref{tab:equilibria} reports the equilibrium solutions obtained by initializing the algorithm from low-price initial pricing vectors. Although several equilibrium prices and the corresponding reinsurance decisions \((A,M)\) exhibit moderate variation across runs, the resulting insolvency levels remain nearly identical, indicating that the optimization drives the insolvency constraint close to its active level. While Table~\ref{tab:equilibria} presents only five representative solutions, the algorithm converged to local equilibrium solutions from all 50 low-price initializations, demonstrating robust convergence under careful initialization. Moreover, all 50 solutions passed the local Nash equilibrium validation test described later, with no profitable deviation identified in their neighborhoods.


\begin{table}[htp]
\centering
\caption{
Equilibrium solutions obtained from five random initializations, where each component of the initial pricing vector is sampled independently from the interval $(0.20,\,0.50)$. The values of \(A\) and \(M\) are reported in units of \(10^8\), while \(\Phi\) denotes the objective value reported in units of \(10^8\). $\alpha$ denotes the insolvency percentage $(\alpha=100 g(\BFx_i,\BFx_{-i}))$.
}
\label{tab:equilibria}

\resizebox{\textwidth}{!}{
\begin{tabular}{c|cccccccc|cc|cc}
\hline
Eq ID
& $\lambda_1$
& $\lambda_2$
& $\lambda_3$
& $\lambda_4$
& $\lambda_5$
& $\lambda_6$
& $\lambda_7$
& $\lambda_8$
& $A$
& $M$
& $\Phi$
& $\alpha$ (\%)
\\
\hline

1
& 0.787 & 0.501 & 0.411 & 1.802 & 1.050 & 0.848 & 1.981 & 0.425
& 2.00 & 9.73 & 8.45 & 0.0463
\\

2
& 0.786 & 0.407 & 0.292 & 1.677 & 0.944 & 0.848 & 1.981 & 0.425
& 1.93 & 10.30 & 7.94 & 0.0463
\\

3
& 0.860 & 0.501 & 0.295 & 1.676 & 0.944 & 0.775 & 1.981 & 0.425
& 1.91 & 10.30 & 7.89 & 0.0463
\\

4
& 0.721 & 0.501 & 0.297 & 1.677 & 1.032 & 0.848 & 1.981 & 0.425
& 1.97 & 10.00 & 8.18 & 0.0463
\\

5
& 0.862 & 0.404 & 0.295 & 1.672 & 0.944 & 0.775 & 1.981 & 0.414
& 1.91 & 10.31 & 7.85 & 0.0463
\\

\hline
\end{tabular}
}
\end{table}

\begin{figure} [htp]
\centering
\subfloat[Convergence trajectories]{%
\resizebox*{6.7cm}{!}{\includegraphics{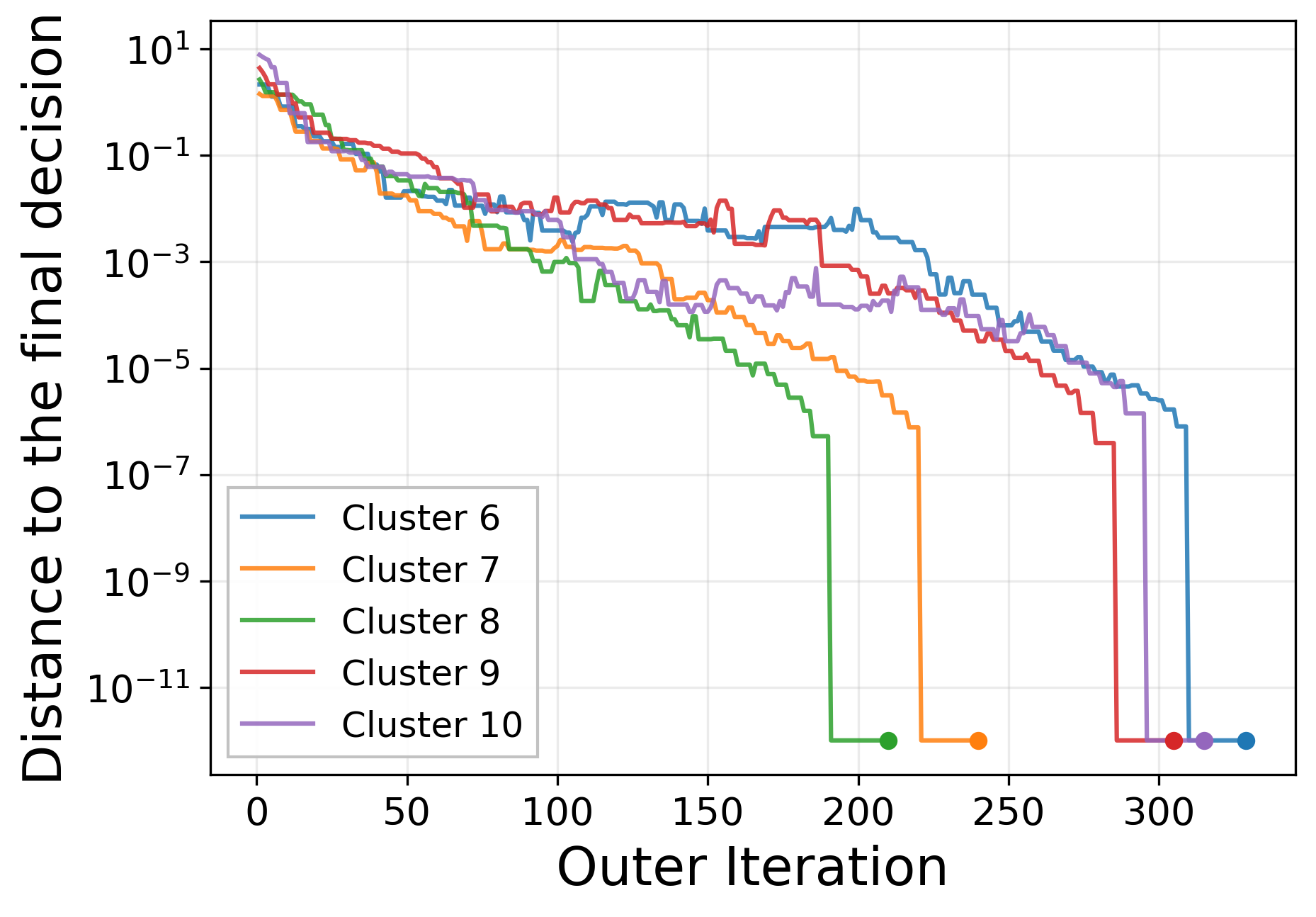}}\label{fig:7a}}
\subfloat[Equilibrium objective values]{%
\resizebox*{6.3cm}{!}{\includegraphics{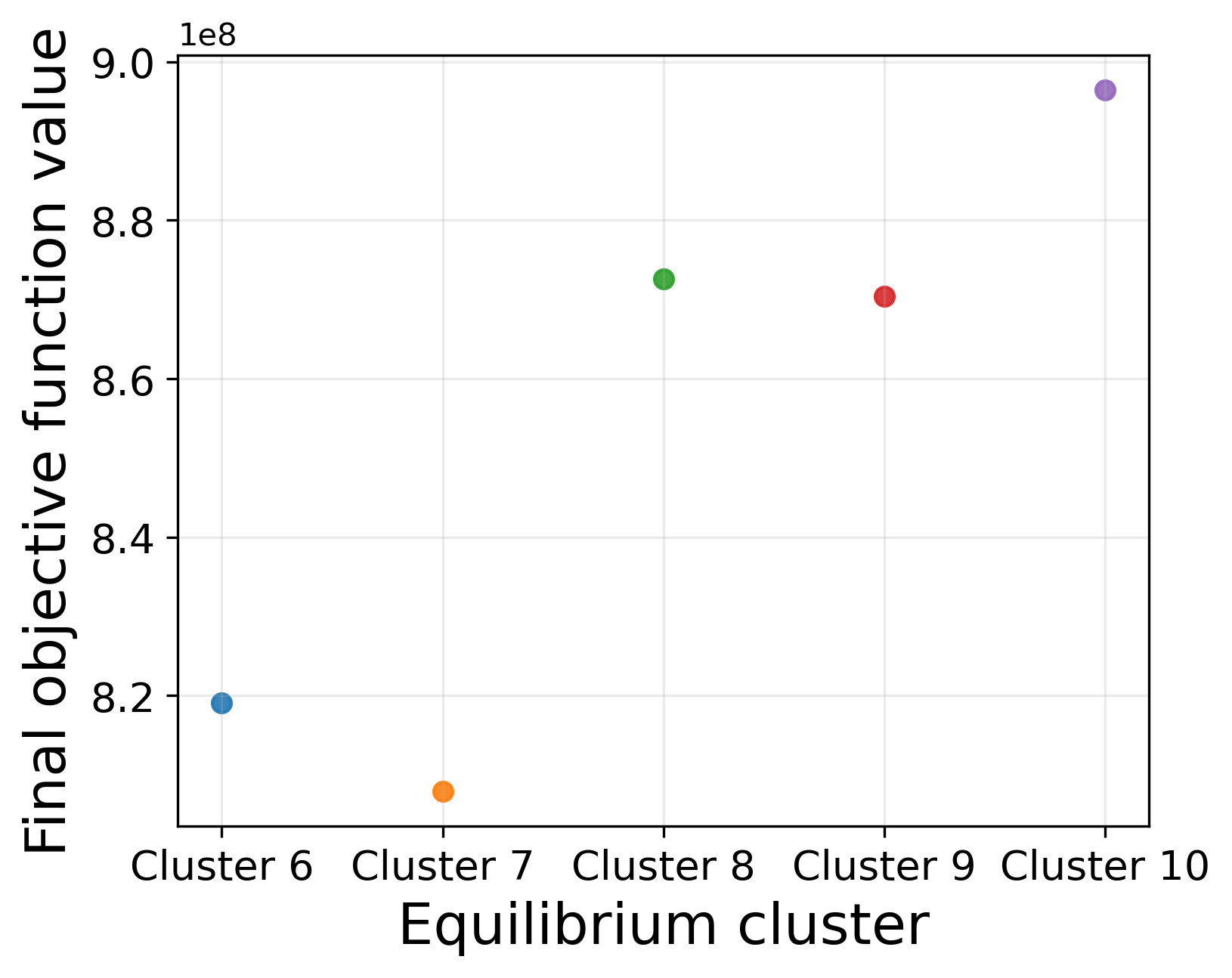}}\label{fig:7b}}
\caption{
Representative equilibrium solutions obtained from different random initializations.
Panel (a) shows the distance to the final solution over the equilibrium iterations, while panel (b) reports the corresponding objective values at the equilibrium solutions.
}
\label{fig:convergence-withoutafford}
\end{figure}

The low-price initialization is intended to provide a reliable mechanism for computing a representative approximate local equilibrium and does not guarantee convergence to a local Nash equilibrium from arbitrary initial pricing vectors. To investigate whether additional equilibrium solutions exist, we next initialize the algorithm from a broader collection of initial pricing vectors. The objective of this experiment is not to exhaustively identify all local Nash equilibria, but rather to demonstrate that multiple local Nash equilibria can arise under different initializations. Specifically, we consider five initialization regions: $(0.1,0.8)$, $(0.6,1.5)$, $(1.2,2.4)$, $(2.0,3.2)$, and $(3.0,4.4)$.
For each region, every component of the initial pricing vector is sampled independently and uniformly from the corresponding interval. Table~\ref{tab:equilibria-multi} reports one representative equilibrium solution obtained from each initialization region, while Figure~\ref{fig:convergence-withoutafford} illustrates the corresponding convergence trajectories and the resulting equilibrium objective values.

\begin{table}[htp]
\centering
\caption{
Representative equilibrium solutions obtained from five initialization regions, where each component of the initial pricing vector is independently sampled from the indicated interval. The notation and units are the same as those in Table~\ref{tab:equilibria}.
}
\label{tab:equilibria-multi}

\resizebox{\textwidth}{!}{
\begin{tabular}{c|cccccccc|cc|cc}
\hline
Eq ID
& $\lambda_1$
& $\lambda_2$
& $\lambda_3$
& $\lambda_4$
& $\lambda_5$
& $\lambda_6$
& $\lambda_7$
& $\lambda_8$
& $A$
& $M$
& $\Phi$
& $\alpha$ (\%)
\\
\hline

6
& 0.718 & 0.501 & 0.293 & 1.746 & 1.015 & 0.848 & 1.981 & 0.425
& 1.97 & 10.00 & 8.19 & 0.0463
\\

7
& 0.830 & 0.501 & 0.297 & 1.676 & 1.021 & 0.779 & 1.981 & 0.425
& 1.95 & 10.01 & 8.08 & 0.0463
\\

8
& 0.803 & 0.497 & 0.411 & 2.139 & 1.050 & 0.825 & 1.981 & 0.414
& 2.05 & 9.43 & 8.73 & 0.0463
\\

9
& 0.678 & 0.407 & 0.295 & 2.229 & 1.015 & 0.775 & 3.221 & 0.414
& 2.07 & 9.72 & 8.70 & 0.0463
\\

10
& 0.660 & 0.404 & 0.274 & 2.194 & 1.148 & 0.775 & 3.221 & 0.421
& 2.11 & 9.44 & 8.96 & 0.0463
\\

\hline
\end{tabular}
}
\end{table}

\begin{remark}
Economically, each equilibrium represents a different joint configuration of premiums and reinsurance decisions. As the overall pricing level increases, the attachment point $A$ generally increases while the exhaustion point $M$ generally decreases, resulting in a narrower reinsurance coverage layer. At the same time, the insurer objective value generally increases, while the insolvency measure remains unchanged. This pattern suggests that the additional premium revenue at the higher-price equilibria allows insurers to retain more catastrophic risk and purchase less reinsurance while continuing to satisfy the insolvency constraint. However, higher premiums may also reduce insurance demand and household participation. Consequently, the equilibria reflect different tradeoffs among insurer financial performance, reinsurance protection, solvency, and insurance accessibility. Quantifying the resulting social welfare implications provides an important direction for future research.
\end{remark}

To assess the computed solutions under the \((\epsilon,r)\) local
Nash equilibrium condition in \eqref{eq:epsilon-r-local-nash}, we consider the problem~\eqref{eq:projected-br}. Although a firm's full decision vector consists of \((\BFlambda_i,A,M)\), many computed solutions lie near the boundary of the insolvency constraint. Direct perturbations of the full decision vector therefore frequently produce infeasible deviations and provide limited information about profitable feasible responses. Instead,
we consider deviations in the pricing vector while reoptimizing the
reinsurance decisions subject to the insolvency constraint.
Specifically, we assess the following projected \((\epsilon,r)\) condition,
\begin{equation}
\begin{aligned}
    f\!\left(
        \BFlambda_i^*,
        \mathcal P(\BFlambda_i^*;\BFx_{-i}^*),
        \BFx_{-i}^*
    \right)
    \geq
    \sup_{\substack{
        \BFlambda_i\in
        \Lambda_i^{\mathrm{feas}}(\BFx_{-i}^*)\\
        \|\BFlambda_i-\BFlambda_i^*\|\leq r
    }}
    f\!\left(
        \BFlambda_i,
        \mathcal P(\BFlambda_i;\BFx_{-i}^*),
        \BFx_{-i}^*
    \right)
    -
    \epsilon.
\end{aligned}
\label{eq:projected-local-nash-test}
\end{equation}
To assess this condition numerically, we sample pricing vectors from the neighborhood in~\eqref{eq:projected-local-nash-test} and use Algorithm~\ref{alg:AM-optimize} to approximate \(\mathcal P(\BFlambda_i;\BFx_{-i}^*)\) for each sample. Thus, each sampled pricing deviation is evaluated after \(A\) and \(M\) have been reoptimized. The projected condition is more stringent than a direct application of~\eqref{eq:epsilon-r-local-nash} to the full decision vector.

For each equilibrium candidate and each radius \(r\in\{0.1,1.0,5.0\}\), we generated 2,000 independent pricing perturbations and reoptimized the reinsurance decisions \(A\) and \(M\) for every perturbed pricing vector. Across the ten candidates, no improving deviation was found among the 20,000 perturbed strategies evaluated at each radius, while demand changes were observed at all three radii. For each candidate-radius pair, zero improvements among the 2,000 samples yield an exact one-sided \(95\%\) Clopper--Pearson upper confidence bound of \(0.15\%\) on the probability that a uniformly sampled pricing perturbation from the corresponding neighborhood improves the objective after \(A\) and \(M\) are reoptimized. One important observation is that no improving deviation was identified even within the largest neighborhood radius. Since the radius \(r=5.0\) covers a substantial portion of the feasible pricing region, the local Nash equilibrium test explores deviations far beyond an infinitesimal neighborhood of the equilibrium solution. Consequently, the identified equilibria appear numerically much closer to best responses than would typically be implied by a local Nash equilibrium condition. Detailed distributions of the objective changes and demand responses are reported in \ref{appdx:s6}.

The average wall-clock time required to compute a single equilibrium solution was \textit{approximately 10 minutes} on a single GPU. Allowing a larger $\tau_{\mathrm{acc}}$ in Algorithm~\ref{alg:lambda-optimize} reduced the average runtime to \textit{approximately 5 minutes}. Despite this relaxed tolerance, no improving deviations were detected when the resulting solutions were evaluated using the same local Nash equilibrium tests; see \ref{appdx:s6}. Given the complexity of the underlying nonconvex equilibrium problem, these computational requirements are modest to support the repeated solution of related game-theoretic models, such as Stackelberg games and sensitivity analyses. Another useful measure of computational effort is the number of objective evaluations required by the equilibrium search procedure. For example, the approach proposed in~\citet{liu2025computing} requires \textit{at least 190,000 objective evaluations} to compute a single equilibrium solution. In contrast, the proposed algorithm typically converges within 300 outer iterations. Since each iteration requires approximately 33 objective evaluations, a representative equilibrium solution can usually be obtained using roughly \textit{9,900 objective evaluations.} This suggests a substantial reduction in computational effort while maintaining high-quality equilibrium solutions.

\subsubsection{With Affordability Cap}

We now present the results under the 5\% affordability cap. Compared with the case without the cap, the equilibrium objective values are substantially lower, indicating a reduction in insurer profitability. However, the cap still permits variation in equilibrium pricing, as illustrated by the distinct solutions obtained from the same initial pricing vector \((0.3,\ldots,0.3)\) in Table~\ref{tab:equilibria-afford}. 
In particular, pricing differences are concentrated in Zones~1--3, suggesting that the scope for alternative equilibrium pricing patterns varies across regions even under a common affordability policy.

\begin{table}[htp]
\centering
\caption{
Equilibrium solutions obtained from the initial pricing vector $(0.3,\ldots,0.3)$ under the 5\% affordability constraint.
The notation is the same as in Table~\ref{tab:equilibria}.
}
\label{tab:equilibria-afford}

\resizebox{\textwidth}{!}{
\begin{tabular}{c|cccccccc|cc|cc}
\hline
Eq ID
& $\lambda_1$ & $\lambda_2$ & $\lambda_3$ & $\lambda_4$
& $\lambda_5$ & $\lambda_6$ & $\lambda_7$ & $\lambda_8$
& $A$ & $M$ & $\Phi$ & $\alpha \,(\%)$ \\
\hline
1
& 0.919 & 0.669 & 0.302 & 0.722 & 0.322 & 0.371 & 0.318 & 0.192
& 1.14 & 59.00 & 3.31 & 0.0408
\\

2
& 1.285 & 0.669 & 0.493 & 0.722 & 0.322 & 0.371 & 0.318 & 0.192
& 1.16 & 59.00 & 3.58 & 0.0408
\\

3
& 0.889 & 0.594 & 0.373 & 0.722 & 0.321 & 0.371 & 0.318 & 0.192
& 1.14 & 59.00 & 3.32 & 0.0408
\\

4
& 0.954 & 0.669 & 0.356 & 0.722 & 0.322 & 0.360 & 0.318 & 0.192
& 1.14 & 59.00 & 3.35 & 0.0408
\\

5
& 0.893 & 0.587 & 0.373 & 0.722 & 0.318 & 0.361 & 0.308 & 0.192
& 1.14 & 59.00 & 3.30 & 0.0408
\\
\hline
\end{tabular}
}
\end{table}

To further investigate the existence of additional equilibrium solutions, we also considered the same five initialization regions as in the case without the affordability cap. Representative solutions are reported in Table~\ref{tab:equilibria-afford-multi}. All ten equilibrium solutions were evaluated using the same local Nash equilibrium verification procedure described previously. Across all ten equilibrium solutions, no improving deviation was found among the 20,000 perturbed strategies evaluated at each radius, providing strong numerical evidence that each identified solution satisfies the local Nash equilibrium condition. As in the baseline case, allowing a modest optimality gap reduces the average wall-clock time under the affordability cap from \textit{approximately 10 minutes} to \textit{approximately 5 minutes}.
Among the ten solutions obtained with this relaxed tolerance, the same local Nash equilibrium tests detected no improving deviations for eight solutions, while the maximum observed relative objective gains for the remaining two were approximately (2.96\%) and (0.89\%); see \ref{appdx:s6}.

\begin{table}[htp]
\centering
\caption{
Equilibrium solutions obtained from five initialization regions under the 5\% affordability constraint, where each component of the initial pricing vector is independently sampled from the corresponding interval. The notation is the same as in Table~\ref{tab:equilibria}.
}
\label{tab:equilibria-afford-multi}

\resizebox{\textwidth}{!}{
\begin{tabular}{c|cccccccc|cc|cc}
\hline
Eq ID
& $\lambda_1$ & $\lambda_2$ & $\lambda_3$ & $\lambda_4$
& $\lambda_5$ & $\lambda_6$ & $\lambda_7$ & $\lambda_8$
& $A$ & $M$ & $\Phi$ & $\alpha \,(\%)$ \\
\hline
6
& 0.921 & 0.576 & 0.350 & 0.722 & 0.316 & 0.360 & 0.318 & 0.192
& 1.14 & 59.00 & 3.29 & 0.0408
\\

7
& 0.869 & 0.503 & 0.314 & 0.722 & 0.303 & 0.371 & 1.027 & 0.192
& 1.14 & 59.00 & 3.27 & 0.0408
\\

8
& 1.285 & 0.653 & 0.413 & 0.722 & 0.320 & 0.371 & 0.272 & 0.192
& 1.15 & 59.00 & 3.51 & 0.0408
\\

9
& 1.049 & 0.669 & 0.356 & 0.722 & 0.303 & 0.371 & 0.318 & 2.662
& 1.14 & 59.00 & 3.35 & 0.0408
\\

10
& 1.285 & 0.530 & 0.356 & 0.722 & 0.304 & 0.371 & 0.318 & 3.449
& 1.15 & 59.00 & 3.39 & 0.0408
\\
\hline
\end{tabular}
}
\end{table}


\section{Conclusion}
In this paper, we studied a simulation-based equilibrium problem for multi-region hurricane insurance markets with multi-period capital dynamics and developed a framework for computing approximate local Nash equilibria under discontinuous demand, reinsurance decisions, and insolvency constraints. Experiments based on a North Carolina hurricane insurance market demonstrate the framework's effectiveness and reveal a complex equilibrium landscape, fundamentally reshaped by interactions between insolvency constraints and reinsurance optimization, with multiple distinct approximate local Nash equilibria. Affordability constraints can induce sharp demand declines over narrow premium ranges and substantially complicate the identification of feasible reinsurance decisions satisfying insolvency constraints. Moreover, the existence and size of the feasible region are highly sensitive to affordability requirements and initial insurer capital, suggesting that market feasibility can vary dramatically across financial and regulatory environments.

Several directions remain for future research. From a methodological perspective, stronger theoretical guarantees for convergence and equilibrium characterization in simulation-based games remain largely unexplored. From an application perspective, the proposed framework can serve as a computational engine for higher-level policy optimization problems, such as Stackelberg formulations involving government interventions, subsidies, or regulatory policies designed to improve market efficiency and social welfare under catastrophic risk.

\section*{Acknowledgements}

Funding: This research was supported by the National Science Foundation under Award No. 2209190. Any opinions, findings, and conclusions expressed in this material are those of the authors and do not necessarily reflect the views of the National Science Foundation.

\appendix


\section{Homeowner Utility and Insurance Demand}
\label{appdx:s1}
The insurance purchase decision of homeowner \(j\) is modeled using the
constant absolute risk aversion utility function
\begin{equation}
    u_j(w)=-\exp(-\theta_j w),
    \label{eq:supp-cara-utility}
\end{equation}
where \(w\) denotes wealth and \(\theta_j>0\) is the homeowner-specific risk
aversion coefficient.

Let \(L_{j,h}^r\) be the loss of homeowner \(j\) under hurricane event type \(h\)
and coverage type \(r\), and let
\(
    B_{j,h}^r=\min\{L_{j,h}^r,2500\}
\)
be the portion retained by the homeowner under the deductible. The maximum
premium that homeowner \(j\) is willing to pay is obtained by equating expected
utility with and without insurance. Under exponential utility, it is
\begin{equation}
    \overline{\pi}_{j}^r
    =
    \frac{1}{\theta_j}
    \log
    \left(
        \frac{
            p_0+
            \displaystyle\sum_{h=1}^{97}
            p_h\exp\!\left(\theta_jL_{j,h}^r\right)
        }{
            p_0+
            \displaystyle\sum_{h=1}^{97}
            p_h\exp\!\left(\theta_jB_{j,h}^r\right)
        }
    \right),
    \label{eq:supp-wtp}
\end{equation}
where
\(
    p_0=1-\sum_{h=1}^{97}p_h
\)
is the probability that none of the 97 hurricane event types occurs during a
given year. Although wealth enters Equation~\eqref{eq:supp-cara-utility}, it
does not appear in Equation~\eqref{eq:supp-wtp} because the common
wealth-dependent term cancels when expected utilities are equated.

The expected insured loss and the corresponding annual premium for
homeowner \(j\) in region \(z\) are given by
\(
    \ell_j^r
    =
    \sum_{h=1}^{97}
    p_h\bigl(L_{j,h}^r-B_{j,h}^r\bigr)
\)
and
\(
    \pi_j^r(\lambda_z)
    =
    (1.35+\lambda_z)\ell_j^r.
\)
Consequently, depending on whether the affordability constraint is imposed,
homeowner \(j\) purchases coverage type \(r\) if
\begin{equation}
    \pi_j^r(\lambda_z)
    \leq
    \begin{cases}
        \overline{\pi}_j^r,
        & \text{without the affordability constraint},\\[1mm]
        \min\!\left\{
            \overline{\pi}_j^r,\,
            0.05H_j
        \right\},
        & \text{with the affordability constraint},
    \end{cases}
    \label{eq:supp-purchase-condition}
\end{equation}
where \(H_j\) denotes the property value of homeowner \(j\). Thus, the exact
purchase decision is represented by an indicator function. Regional demand,
premium revenue, retained deductible, and hurricane-loss exposure are obtained
by aggregating the corresponding homeowner-level quantities within each region.

\FloatBarrier

\section{Capital Dynamics and Cash-Flow Components}
\label{appdx:s2}

Let \(i\) denote an insurer, \(y=1,\ldots,T\) the year, and
\(s\in\mathcal{S}\) a hurricane scenario path. For strategy profile
\(\BFx\), let
\(\widetilde C_{i,y}^s(\BFx)\) and \(C_{i,y}^s(\BFx)\) denote,
respectively, insurer \(i\)'s capital at the end of year \(y\) before
and after any external capital injection. Before an injection, capital
evolves according to
\begin{equation}
    \widetilde C_{i,y}^{s}(\BFx)
    =
    C_{i,y-1}^{s}(\BFx)
    +
    \beta^{y-1}
    \left[
        \Pi_i(\BFx)
        +
        D_{i,y}^{s}(\BFx)
        -
        L_{i,y}^{s}(\BFx)
        +
        R_{i,y}^{s}(\BFx)
    \right],
    \label{eq:supp-capital-dynamics}
\end{equation}
where \(\beta=0.95\) is the annual discount factor. Initial capital is set to
\(
    C_{i,0}^{s}(\BFx)
    =
    3P_i(\BFx),
\)
where \(P_i(\BFx)\) is insurer \(i\)'s gross annual premium revenue.
If \(\widetilde C_{i,y}^{s}(\BFx)\leq0\), the incumbent insurer exits
and a replacement insurer enters under the same pricing and reinsurance strategy. The associated external capital injection is
\[
    I_{i,y}^{s}(\BFx)
    =
    3\bigl(\Pi_i(\BFx)\bigr)_+
    \mathbf{1}
    \left\{
        \widetilde C_{i,y}^{s}(\BFx)\leq0
    \right\},
\]
where \((a)_+:=\max\{a,0\}\), and the post-injection capital is
\(
    C_{i,y}^{s}(\BFx)
    =
    \widetilde C_{i,y}^{s}(\BFx)
    +
    I_{i,y}^{s}(\BFx).
\)

\paragraph{Annual premium income.}
Let
\(\ell_{i,z}(\BFx)\)
denote insurer \(i\)'s expected insured loss exposure in region \(z\),
as determined by the household demand model under strategy profile
\(\BFx\). The gross annual premium revenue is
\begin{equation}
    P_i(\BFx)
    =
    \sum_{z=1}^{8}
    (1.35+\lambda_{i,z})
    \ell_{i,z}(\BFx).
    \label{eq:supp-gross-premium}
\end{equation}
After deducting non-claim operating costs equal to \(35\%\) of expected
insured loss, annual premium income net of operating costs is
\begin{align}
    \Pi_i(\BFx)
    =
    P_i(\BFx)
    -
    0.35
    \sum_{z=1}^{8}
    \ell_{i,z}(\BFx)
    =
    \sum_{z=1}^{8}
    (1+\lambda_{i,z})
    \ell_{i,z}(\BFx).
    \label{eq:supp-net-premium}
\end{align}

\paragraph{Realized losses and deductibles.}
Let \(\mathcal{H}_y^s\) denote the set of hurricane events occurring in
year \(y\) under scenario path \(s\). For homeowner \(j\), coverage type
\(r\in\{\mathrm{wind},\mathrm{flood}\}\), and hurricane event type \(h\),
let \(L_{j,h}^r\) denote the property loss before application of the
deductible. The portion of the loss retained by the homeowner under the
\$2,500 deductible is
\(
    B_{j,h}^r
    =
    \min
    \{
        L_{j,h}^r,\,
        2500
    \}.
\)
For a pricing vector \(\boldsymbol{\lambda}_i\), the household demand
model determines whether each homeowner purchases wind and flood
coverage. Let
\(
    a_{j}^{r}(\boldsymbol{\lambda}_i)\in\{0,1\}
\)
denote the corresponding purchase indicator for coverage type \(r\).
The indicator equals one when homeowner \(j\) satisfies the purchase conditions under \(\boldsymbol{\lambda}_i\) and zero otherwise. 
For region \(z\), let \(\mathcal J_z\) denote the set of homeowners in the region. Then the aggregate event-level loss and deductible amounts associated with households satisfying the purchase conditions under
\(\boldsymbol{\lambda}_i\) are given by
\begin{align}
    L_{z,h}(\boldsymbol{\lambda}_i)
    &=
    \sum_{r\in\{\mathrm{wind},\mathrm{flood}\}}
    \sum_{j\in\mathcal J_z}
    a_j^r(\boldsymbol{\lambda}_i)L_{j,h}^r,
    \label{eq:supp-zone-event-loss}\\
    D_{z,h}(\boldsymbol{\lambda}_i)
    &=
    \sum_{r\in\{\mathrm{wind},\mathrm{flood}\}}
    \sum_{j\in\mathcal J_z}
    a_j^r(\boldsymbol{\lambda}_i)B_{j,h}^r.
    \label{eq:supp-zone-event-deductible}
\end{align}

Competition determines the portion of the resulting regional portfolio
served by insurer \(i\). Let \(d_z(\boldsymbol{\lambda}_i)\) denote the
total insured loss exposure generated by the household demand model in
region \(z\), and let \(d_{-i,z}(\BFx_{-i})\) denote the exposure served
by competing insurers. Insurer \(i\)'s insured exposure in region \(z\) then becomes
\(
    \ell_{i,z}(\BFx)
    =
    \left[
        d_z(\boldsymbol{\lambda}_i)
        -
        d_{-i,z}(\BFx_{-i})
    \right]_+.
\)
When \(d_z(\boldsymbol{\lambda}_i)>0\), define the corresponding
portfolio share as
\(
    \alpha_{i,z}(\BFx)
    =
    \ell_{i,z}(\BFx)/
         {d_z(\boldsymbol{\lambda}_i)},
\)
with \(\alpha_{i,z}(\BFx)=0\) when
\(d_z(\boldsymbol{\lambda}_i)=0\).
The same regional portfolio share is applied to event-level losses and
deductibles. The annual quantities under scenario path \(s\) are then given by
\begin{align*}
    L_{i,y}^s(\BFx)
    &=
    \sum_{h\in\mathcal H_y^s}
    \sum_{z=1}^{8}
    \alpha_{i,z}(\BFx)
    L_{z,h}(\boldsymbol{\lambda}_i),\\
    D_{i,y}^s(\BFx)
    &=
    \sum_{h\in\mathcal H_y^s}
    \sum_{z=1}^{8}
    \alpha_{i,z}(\BFx)
    D_{z,h}(\boldsymbol{\lambda}_i).
\end{align*}
Consequently, insurer \(i\)'s realized claim payment before reinsurance is
\(L_{i,y}^s(\BFx)-D_{i,y}^s(\BFx)\).

\paragraph{Net reinsurance cash flow.}
Let \(L_{i,h}(\BFx)\) denote insurer \(i\)'s aggregate insured loss
under hurricane event \(h\) for strategy profile \(\BFx\). The gross
payout from insurer \(i\)'s reinsurance layer \([A_i,M_i]\) is
\[
    q_{i,h}(\BFx)
    =
    \min\left\{
        \max\left\{
            L_{i,h}(\BFx)-A_i,\,
            0
        \right\},
        M_i-A_i
    \right\}.
\]
Because the contract is applied separately to each hurricane event, the
gross reinsurance payout in year \(y\) under scenario path \(s\) is
\(
    Q_{i,y}^s(\BFx)
    =
    \sum_{h\in\mathcal{H}_y^s}
    q_{i,h}(\BFx).
\)

Let \(p_h\) denote the annual occurrence probability of hurricane event
type \(h\). For reinsurance pricing, define the probability-weighted
payout and a payout dispersion measure as
\begin{align*}
    \overline{q}_i(\BFx)
    &=
    \sum_{h=1}^{97}
    p_h q_{i,h}(\BFx),\\
    \sigma_i(\BFx)
    &=
    \left[
        \sum_{h=1}^{97}
        p_h
        \left(
            q_{i,h}(\BFx)
            -
            \overline{q}_i(\BFx)
        \right)^2
    \right]^{1/2}.
\end{align*}
Here, \(\overline{q}_i\) is the probability-weighted sum of the
event-level payouts, and \(\sigma_i\) measures their dispersion
around \(\overline{q}_i\) over the 97 hurricane event types using the
same event weights. These two quantities enter the reinsurance pricing
rule through the fixed annual reinsurance charge \(b_i\) and the
variable reinsurance charge rate \(v_i\):
\begin{align*}
    b_i(\BFx)
    &=
    (1+\phi)\beta
    \overline{q}_i(\BFx)
    +
    g\beta
    \sigma_i(\BFx),\\
    v_i(\BFx)
    &=
    \frac{
        \beta\overline{q}_i(\BFx)
    }{
        M_i-A_i
    },
\end{align*}
where \(\phi=g=0.1\). The fixed charge \(b_i\) is paid each year
regardless of the realized reinsurance payout. Its first component,
\((1+\phi)\beta\overline{q}_i\), applies a loading of \(\phi\) to the
scaled probability-weighted payout, while the second component,
\(g\beta\sigma_i\), adds a risk loading proportional to the payout
dispersion measure. Thus, ``fixed'' refers to independence from the
realized payout, rather than independence from the pricing and
reinsurance decisions. The coefficient \(v_i\) is a variable-charge
rate obtained by dividing the scaled probability-weighted payout by
the width of the reinsurance layer, \(M_i-A_i\). The corresponding
charge is calculated by multiplying \(v_i\) by the realized annual
payout \(Q_{i,y}^s\), and therefore varies in direct proportion to that
payout.

Consequently, the net reinsurance cash flow in
\eqref{eq:supp-capital-dynamics} is
\begin{align*}
    R_{i,y}^s(\BFx)
    =
    \beta Q_{i,y}^s(\BFx)
    -
    \left[
        b_i(\BFx)
        +
        v_i(\BFx)
        Q_{i,y}^s(\BFx)
    \right],
\end{align*}
where the first term is the realized payout scaled by \(\beta\), and
the two terms in brackets represent the fixed annual charge and the
variable charge proportional to the realized payout, respectively.

\section{Finite-Difference Approximation}
\label{appdx:s3}

This section describes the finite-difference (FD) approximations used in the
price-optimization procedure of the main manuscript. At iteration \(k\), we
generate a \(2d+1\) FD stencil within the trust-region neighborhood of
\(\BFlambda_k\). Let
\(
    \{v_i\}_{i=1}^d
\)
be an orthonormal collection of search directions. Depending on the algorithmic
setting, these directions are either randomized or chosen as the coordinate
basis. For each \(v_i\), we evaluate
\(
    \BFlambda_k^{(i,+)}
    =
    \BFlambda_k+h_i^+v_i,
    \BFlambda_k^{(i,-)}
    =
    \BFlambda_k-h_i^-v_i,
\)
where \(h_i^+>0\) and \(h_i^->0\) are feasible step lengths determined by the
trust-region radius and the variable bounds.
When \(h_i^+=h_i^-\) and the directions form an orthonormal basis, this is the
standard symmetric \(2d+1\) interpolation geometry used in underdetermined
quadratic derivative-free models. \citet{ragonneau2024optimal} establish an
optimal well-poisedness property for the corresponding standard interpolation
set. Near a bound, the use of unequal feasible step lengths preserves the same
directional stencil structure, although the optimality result for the symmetric
geometry does not directly apply.
Define
\(
    f_k=\bar{\Phi}(\BFlambda_k),
    f_k^{(i,+)}=\bar{\Phi}(\BFlambda_k^{(i,+)}),
\) and 
\(
    f_k^{(i,-)}=\bar{\Phi}(\BFlambda_k^{(i,-)}).
\)
The first derivative of the quadratic interpolant along \(v_i\) at
\(\BFlambda_k\) is
\begin{equation}
    g_i^{\mathrm{dir}}
    =
    \frac{
        (h_i^-)^2\bigl(f_k^{(i,+)}-f_k\bigr)
        +(h_i^+)^2\bigl(f_k-f_k^{(i,-)}\bigr)
    }{
        h_i^+h_i^-(h_i^++h_i^-)
    }.
    \label{eq:supp-directional-gradient}
\end{equation}
The corresponding second derivative is
\begin{equation}
    h_i^{\mathrm{dir}}
    =
    \frac{
        2\left[
            h_i^-\bigl(f_k^{(i,+)}-f_k\bigr)
            +h_i^+\bigl(f_k^{(i,-)}-f_k\bigr)
        \right]
    }{
        h_i^+h_i^-(h_i^++h_i^-)
    }.
    \label{eq:supp-directional-curvature}
\end{equation}
In particular, when \(h_i^+=h_i^-\),
\eqref{eq:supp-directional-gradient} and
\eqref{eq:supp-directional-curvature} reduce to the usual centered
finite-difference formulas.

Let
\(
    V=[v_1,\ldots,v_d]
\),
\(
    g^{\mathrm{dir}}
    =(g_1^{\mathrm{dir}},\ldots,g_d^{\mathrm{dir}})^{\mathsf T}
\), and
\(
    H^{\mathrm{dir}}
    =\operatorname{diag}(h_1^{\mathrm{dir}},\ldots,h_d^{\mathrm{dir}})
\).
The directional approximations can be represented in the original coordinate
system as
\(
    \widehat g_k=Vg^{\mathrm{dir}} 
\)
and
\(
    \widehat H_k=VH^{\mathrm{dir}}V^{\mathsf T}.
\)
Thus, the model is diagonal in the orthonormal directional basis, while \(\widehat H_k\) need not be diagonal in the original coordinates. The stencil directions are also reused as candidate direct-search directions. Randomizing the orthonormal basis broadens the directional exploration of the search space. Related analyses of discontinuous direct-search and nonsmooth trust-region methods use increasingly rich or asymptotically dense direction sets to obtain generalized stationarity results~\citep{vicente2012discontinuous,liuzzi2019trust}.
\FloatBarrier

\section{Additional Results under a Larger $\tau_{\text{acc}}$}
\label{appdx:s4}

This section reports additional numerical results obtained using a looser tolerance than in the primary experiments. Specifically, we set \(\tau_{\mathrm{acc}}=10^{-4}\) in Algorithm~2 of the main manuscript. Tables~\ref{tab:supp-equil} and \ref{tab:supp-equilibria-afford} report the resulting solutions without and with the \(5\%\) affordability constraint, respectively. Without an affordability cap, no improving unilateral deviation was identified in the neighborhood tests. Under the 5\% affordability cap, tests using 2,000 perturbations per solution at each radius identified improving deviations for two of the ten reported solutions, both at \(r=0.1\), with maximum observed relative objective gains of approximately \(2.96\%\) and \(0.89\%\). No improving deviations were detected for the remaining eight solutions or at the larger tested radii. Further details of the testing procedure and the distributions of observed relative objective changes are provided in Section~\ref{appdx:s6}.

\begin{table}[htp]
\centering
\caption{Equilibrium solutions obtained from different random initializations
without the affordability constraint, using
\(\tau_{\mathrm{acc}}=10^{-4}\). The variables and scaling follow the notation
of the main manuscript.}
\label{tab:supp-equil}
\resizebox{\textwidth}{!}{%
\begin{tabular}{c|cccccccc|cc|cc}
\hline
Eq ID
& $\lambda_1$ & $\lambda_2$ & $\lambda_3$ & $\lambda_4$
& $\lambda_5$ & $\lambda_6$ & $\lambda_7$ & $\lambda_8$
& $A$ & $M$ & $\Phi$ & $\alpha\,(\%)$ \\
\hline
1  & 0.729 & 0.502 & 0.349 & 1.672 & 1.015 & 0.847 & 1.979 & 0.424 & 1.96 & 10.00 & 8.17 & 0.0463 \\
2  & 0.870 & 0.507 & 0.293 & 1.825 & 1.020 & 0.848 & 1.980 & 0.421 & 1.99 &  9.71 & 8.34 & 0.0463 \\
3  & 0.723 & 0.501 & 0.411 & 1.677 & 1.015 & 0.847 & 1.980 & 0.412 & 1.97 & 10.00 & 8.21 & 0.0463 \\
4  & 0.743 & 0.407 & 0.411 & 1.690 & 1.015 & 0.848 & 1.980 & 0.518 & 1.97 & 10.00 & 8.21 & 0.0463 \\
5  & 0.864 & 0.669 & 0.411 & 2.229 & 1.050 & 0.848 & 1.978 & 0.424 & 2.08 &  9.14 & 8.93 & 0.0463 \\
6  & 0.731 & 0.435 & 0.349 & 1.680 & 1.031 & 0.847 & 1.979 & 0.424 & 1.97 & 10.00 & 8.19 & 0.0463 \\
7  & 0.719 & 0.407 & 0.295 & 1.672 & 1.015 & 0.775 & 1.980 & 0.412 & 1.94 & 10.30 & 7.96 & 0.0463 \\
8  & 0.702 & 0.412 & 0.352 & 2.228 & 1.032 & 0.826 & 3.213 & 0.422 & 2.09 &  9.45 & 8.89 & 0.0463 \\
9  & 0.716 & 0.501 & 0.349 & 1.656 & 1.009 & 0.904 & 1.905 & 0.384 & 1.97 & 10.00 & 8.20 & 0.0463 \\
10 & 0.780 & 0.403 & 0.295 & 1.676 & 1.148 & 0.775 & 3.220 & 0.424 & 2.04 &  9.72 & 8.56 & 0.0463 \\
\hline
\end{tabular}%
}
\end{table}

\begin{table}[htp]
\centering
\caption{Equilibrium solutions obtained from different random initializations
under the \(5\%\) affordability constraint, using
\(\tau_{\mathrm{acc}}=10^{-4}\). The variables and scaling follow the notation
of the main manuscript.}
\label{tab:supp-equilibria-afford}
\resizebox{\textwidth}{!}{%
\begin{tabular}{c|cccccccc|cc|cc}
\hline
Eq ID
& $\lambda_1$ & $\lambda_2$ & $\lambda_3$ & $\lambda_4$
& $\lambda_5$ & $\lambda_6$ & $\lambda_7$ & $\lambda_8$
& $A$ & $M$ & $\Phi$ & $\alpha\,(\%)$ \\
\hline
1  & 0.880 & 0.507 & 0.301 & 0.722 & 0.302 & 0.370 & 0.307 & 0.181 & 1.14 & 59.00 & 3.19 & 0.0408 \\
2  & 0.884 & 0.536 & 0.305 & 0.722 & 0.322 & 0.371 & 0.315 & 0.053 & 1.14 & 59.00 & 3.25 & 0.0408 \\
3  & 0.965 & 0.668 & 0.441 & 0.722 & 0.318 & 0.371 & 0.309 & 0.043 & 1.15 & 59.00 & 3.42 & 0.0408 \\
4  & 0.847 & 0.501 & 0.314 & 0.722 & 0.302 & 0.371 & 1.027 & 0.049 & 1.14 & 59.00 & 3.25 & 0.0408 \\
5  & 0.297 & 0.173 & 0.179 & 0.707 & 0.293 & 0.356 & 0.931 & 0.054 & 1.12 & 59.00 & 2.60 & 0.0408 \\
6  & 0.556 & 0.299 & 0.268 & 0.716 & 0.300 & 0.371 & 1.026 & 0.035 & 1.13 & 59.00 & 2.98 & 0.0408 \\
7  & 1.285 & 0.414 & 0.411 & 0.722 & 0.321 & 0.371 & 0.308 & 0.016 & 1.15 & 59.00 & 3.43 & 0.0408 \\
8  & 1.044 & 0.668 & 0.350 & 0.722 & 0.304 & 0.371 & 0.317 & 2.660 & 1.14 & 59.00 & 3.35 & 0.0408 \\
9  & 1.287 & 0.669 & 0.487 & 0.719 & 0.336 & 0.371 & 4.396 & 3.979 & 1.14 & 59.00 & 3.57 & 0.0408 \\
10 & 1.270 & 0.507 & 0.410 & 0.722 & 0.303 & 0.371 & 0.316 & 3.452 & 1.15 & 59.00 & 3.41 & 0.0408 \\
\hline
\end{tabular}%
}
\end{table}

\FloatBarrier

\section{Sensitivity to the Smoothing Parameter}
\label{appdx:s5}

The parameter \(\mu\) smooths the discontinuous demand function, producing a
continuously differentiable approximation for use within the local optimization
steps. Because \(\mu\) is user specified, we examine whether the computed
equilibrium solutions are sensitive to its value while holding all other model
parameters fixed.

\begin{figure}[htp]
\centering
\subfloat[Full view]{%
\resizebox*{6.7cm}{!}{%
\includegraphics{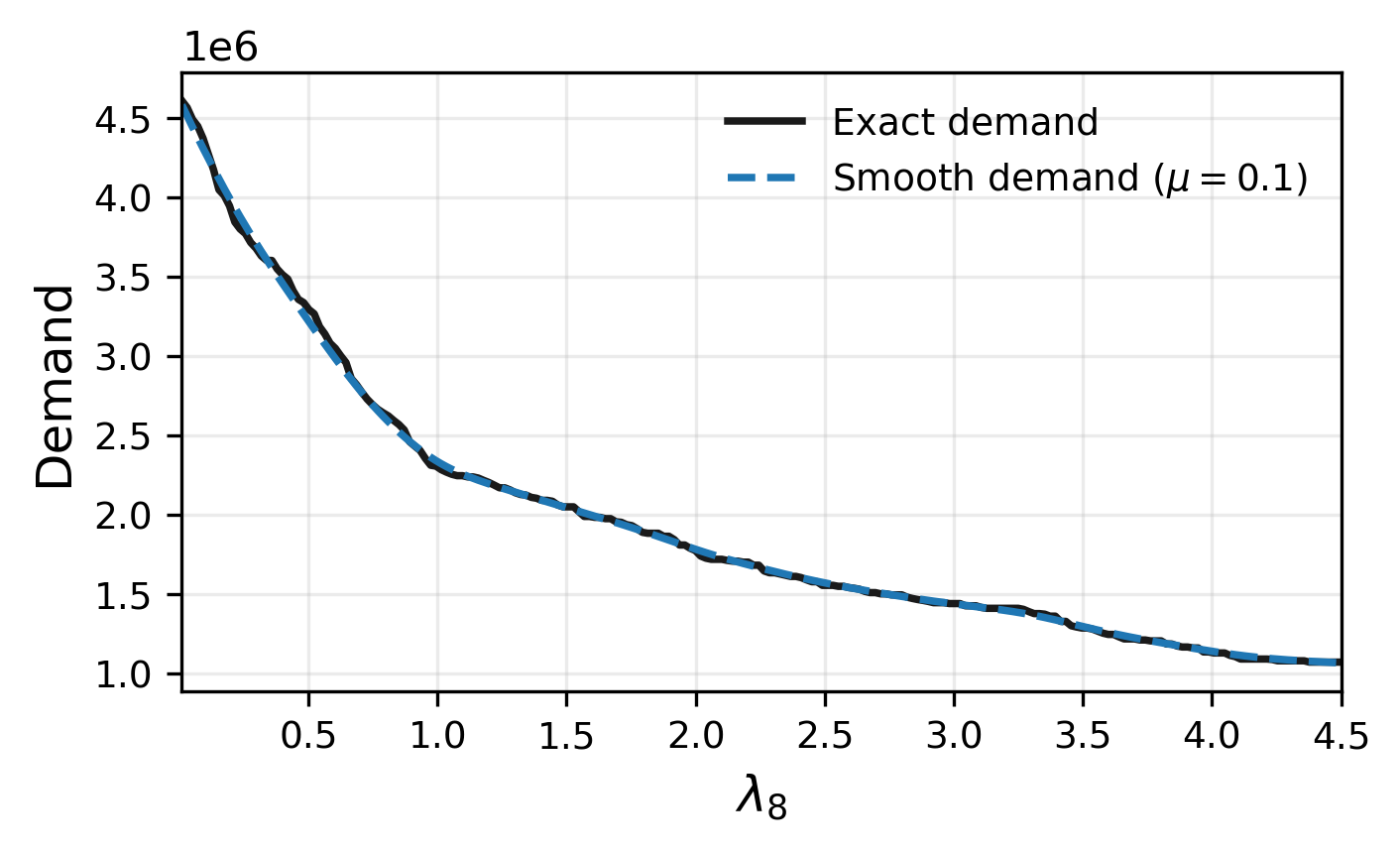}}%
\label{fig:supp-demand-full}}%
\hspace{2pt}
\subfloat[Zoomed-in view]{%
\resizebox*{6.7cm}{!}{%
\includegraphics{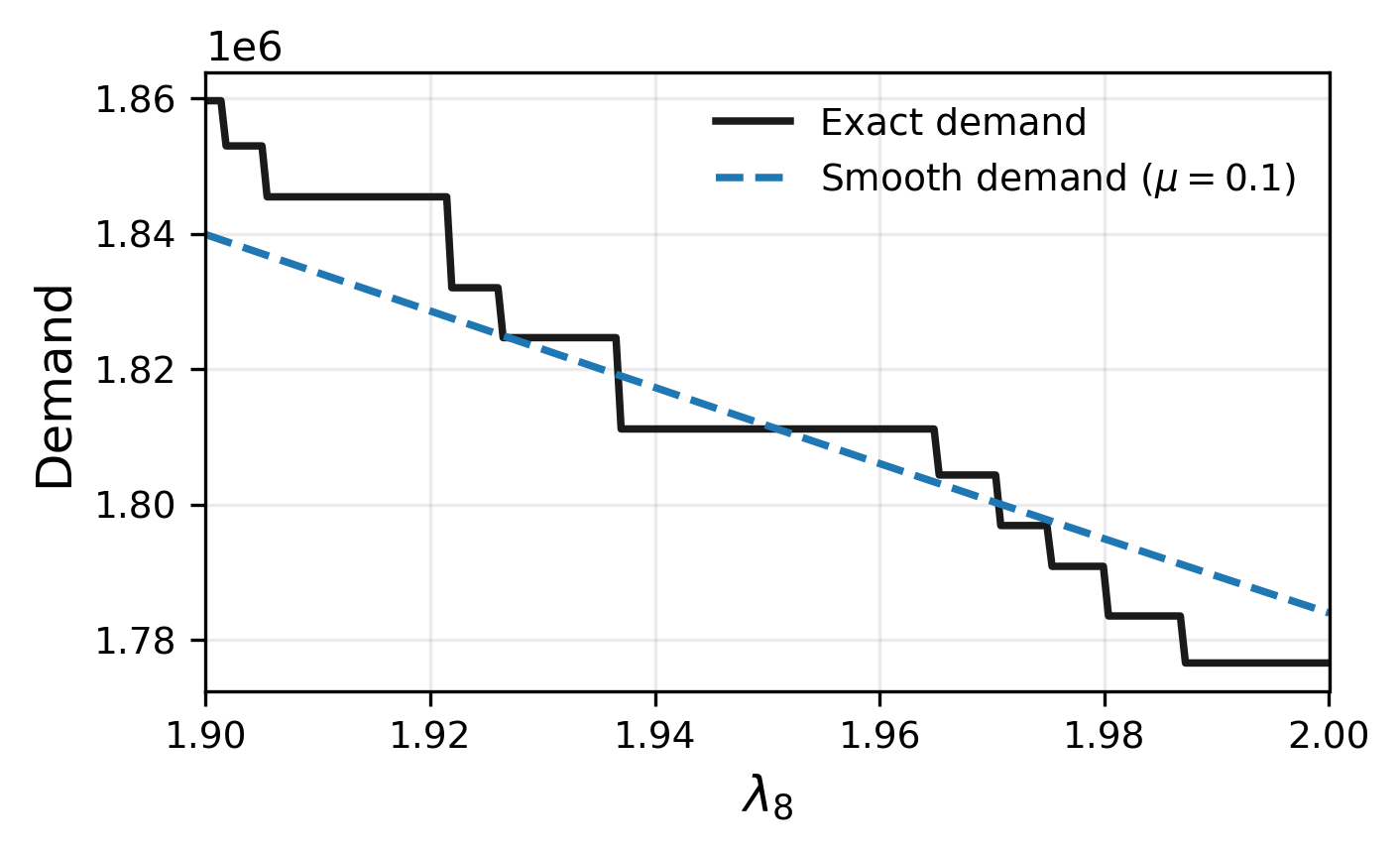}}%
\label{fig:supp-demand-zoom}}
\caption{
Comparison of the original discontinuous demand function and its smoothed approximation with \(\mu=0.1\).
The full view illustrates their overall agreement, while the zoomed-in view highlights smoothing near the transition regions.
}
\label{fig:supp-demand-functions}
\end{figure}

Figure~\ref{fig:supp-demand-functions} compares the original discontinuous
demand function with its smoothed approximation for \(\mu=0.1\). The full view
shows close overall agreement, while the zoomed view shows that each jump is
replaced by an approximately linear transition over a small neighborhood.

We first compute equilibrium solutions using only the smoothed demand function
with \(\mu=0.1\), without reevaluating the candidates under the original
discontinuous demand model. Compared with the baseline equilibria reported in
the main manuscript, the resulting equilibria differ noticeably, indicating
that \(\mu=0.1\) can introduce non-negligible approximation bias when the
smoothed model is treated as the final demand model. However, the five runs in
Table~\ref{tab:supp-equilibria-smoothed-only} form only two practically distinct
solution groups, which differ primarily in the price for Zone~7. In contrast,
the original discontinuous demand function produces a substantially more
diverse collection of equilibrium candidates.

Replacing jumps by smooth transitions changes the local objective landscape and
may remove local optima. These results suggest that larger values of \(\mu\) can
reduce the number of attainable solutions at the cost of greater approximation
bias. Smaller values of \(\mu\) more closely approximate the discontinuous
demand function but preserve a richer set of local solutions.

\begin{remark}[Demand-model selection]
The objective of the proposed framework is to compute equilibria for a
user-specified, data-driven demand model rather than to advocate a particular
demand representation. The method can be applied when each firm's
better-response subproblem can be solved under the selected demand model,
including the original discontinuous demand function, a smoothed approximation,
or a simpler surrogate such as a linear demand model. Equilibria obtained from
a smoothed representation remain equilibria of the approximating model and need
not coincide with equilibria of the original discontinuous model. The demand
representation should therefore be selected according to the intended economic
interpretation and the acceptable approximation error.
\end{remark}

\begin{table}[htp]
\centering
\caption{Equilibrium solutions obtained using only the smoothed demand function
with \(\mu=0.1\). The variables and scaling follow the notation of the main
manuscript.}
\label{tab:supp-equilibria-smoothed-only}
\resizebox{\textwidth}{!}{%
\begin{tabular}{c|cccccccc|cc|cc}
\hline
Eq ID
& $\lambda_1$ & $\lambda_2$ & $\lambda_3$ & $\lambda_4$
& $\lambda_5$ & $\lambda_6$ & $\lambda_7$ & $\lambda_8$
& $A$ & $M$ & $\Phi$ & $\alpha\,(\%)$ \\
\hline
1 & 0.725 & 0.494 & 0.385 & 2.072 & 1.020 & 0.843 & 1.963 & 0.331 & 2.03 & 9.73 & 8.55 & 0.0463 \\
2 & 0.725 & 0.494 & 0.385 & 2.072 & 1.020 & 0.843 & 1.964 & 0.331 & 2.03 & 9.73 & 8.55 & 0.0463 \\
3 & 0.722 & 0.494 & 0.385 & 2.073 & 1.020 & 0.843 & 1.965 & 0.332 & 2.03 & 9.73 & 8.55 & 0.0463 \\
4 & 0.778 & 0.549 & 0.386 & 2.087 & 1.016 & 0.813 & 3.135 & 0.329 & 2.06 & 9.43 & 8.77 & 0.0463 \\
5 & 0.742 & 0.507 & 0.394 & 2.031 & 1.031 & 0.854 & 3.137 & 0.338 & 2.05 & 9.43 & 8.80 & 0.0463 \\
\hline
\end{tabular}%
}
\end{table}

Finally, Table~\ref{tab:supp-equil-mu-small} reports results for
\(\mu=0.01\). Relative to the \(\mu=0.1\) results in
Table~\ref{tab:supp-equil}, the algorithm returns a more diverse collection of
solutions. The solutions for \(\mu=0.01\) span substantially different pricing
and reinsurance decisions and exhibit a wider range of objective values. This
additional diversity is consistent with the smaller smoothing neighborhood
preserving more of the local structure of the original discontinuous demand
model.

\begin{table}[htp]
\centering
\caption{Equilibrium solutions obtained using \(\mu=0.01\) from different
random initializations. The variables and scaling follow the notation of the
main manuscript.}
\label{tab:supp-equil-mu-small}
\resizebox{\textwidth}{!}{%
\begin{tabular}{c|cccccccc|cc|cc}
\hline
Eq ID
& $\lambda_1$ & $\lambda_2$ & $\lambda_3$ & $\lambda_4$
& $\lambda_5$ & $\lambda_6$ & $\lambda_7$ & $\lambda_8$
& $A$ & $M$ & $\Phi$ & $\alpha\,(\%)$ \\
\hline
1  & 0.785 & 0.401 & 0.292 & 1.690 & 0.943 & 0.790 & 1.905 & 0.420 & 1.92 & 10.30 &  7.84 & 0.0463 \\
2  & 0.859 & 0.407 & 0.411 & 1.802 & 0.940 & 0.994 & 1.979 & 0.485 & 2.00 &  9.73 &  8.42 & 0.0463 \\
3  & 0.892 & 0.769 & 0.411 & 1.680 & 1.015 & 0.790 & 1.907 & 0.412 & 1.96 &  9.72 &  8.26 & 0.0463 \\
4  & 0.918 & 0.528 & 0.297 & 1.672 & 1.009 & 0.790 & 1.370 & 0.424 & 1.92 & 10.02 &  7.96 & 0.0463 \\
5  & 1.726 & 2.285 & 0.915 & 1.690 & 1.241 & 0.975 & 1.979 & 0.523 & 2.26 &  8.27 & 10.00 & 0.0463 \\
6  & 0.625 & 0.668 & 0.420 & 1.677 & 0.942 & 1.167 & 1.906 & 0.413 & 2.02 &  9.73 &  8.54 & 0.0463 \\
7  & 1.315 & 0.857 & 0.518 & 2.383 & 1.032 & 0.994 & 2.016 & 0.415 & 2.16 &  8.57 &  9.54 & 0.0463 \\
8  & 2.328 & 1.990 & 0.925 & 2.498 & 1.979 & 1.043 & 3.220 & 0.837 & 4.51 &  6.75 & 13.04 & 0.0456 \\
9  & 3.814 & 2.469 & 1.426 & 3.518 & 2.112 & 1.191 & 3.225 & 0.714 & 5.95 &  7.10 & 14.79 & 0.0456 \\
10 & 3.751 & 3.168 & 1.122 & 3.704 & 1.500 & 1.274 & 3.216 & 1.192 & 5.71 &  7.35 & 13.95 & 0.0456 \\
\hline
\end{tabular}%
}
\end{table}

\FloatBarrier

\section{Distributional Diagnostics for the Local Nash Equilibrium Tests}
\label{appdx:s6}

This section supplements the local Nash equilibrium tests by examining the distributions of objective changes generated by sampled unilateral pricing deviations. We report results separately for the settings without and with the affordability cap, using both the acceptance threshold employed in the main experiments and the larger threshold \(\tau_{\mathrm{acc}}=10^{-4}\). For each perturbed pricing vector, the reinsurance decisions \(A\) and \(M\) were reoptimized using the same procedure employed in the equilibrium computation, while the competitors' decisions were held fixed.

Let \(\BFlambda_{i,c}^{*}\) and \(\BFx_{-i,c}^{*}\) denote insurer \(i\)'s pricing vector and its competitors' decisions at reported solution \(c\), respectively. For a perturbed pricing vector \(\widetilde{\BFlambda}_{i,c,m}^{\,r}\), define the relative objective change as
\begin{equation}
    \delta_{c,m}^{\,r}
    =
    \frac{
        \begin{aligned}
        \bar\Phi\!\left(
            \widetilde{\BFlambda}_{i,c,m}^{\,r},
            \widetilde{\mathcal P}(
                \widetilde{\BFlambda}_{i,c,m}^{\,r};
                \BFx_{-i,c}^{*}
            ),
            \BFx_{-i,c}^{*}
        \right)
        -
        \bar\Phi\!\left(
            \BFlambda_{i,c}^{*},
            \widetilde{\mathcal P}(
                \BFlambda_{i,c}^{*};
                \BFx_{-i,c}^{*}
            ),
            \BFx_{-i,c}^{*}
        \right)
        \end{aligned}
    }{
        \max\left\{
            1,\,
            \left|
                \bar\Phi\!\left(
                    \BFlambda_{i,c}^{*},
                    \widetilde{\mathcal P}(
                        \BFlambda_{i,c}^{*};
                        \BFx_{-i,c}^{*}
                    ),
                    \BFx_{-i,c}^{*}
                \right)
            \right|
        \right\}
    }.
    \label{eq:supp-relative-objective-change}
\end{equation}
Here, \(\bar\Phi\) denotes the expected penalty-augmented objective,
and \(\widetilde{\mathcal P}\) approximates the reinsurance
optimization operator \(\mathcal P\) using Algorithm~1
of the main manuscript, returning approximately optimal
reinsurance parameters for the given pricing vector and
fixed competitors' decisions.
Hence, a positive value of \(\delta_{c,m}^{\,r}\) indicates an
improving unilateral deviation after reoptimizing reinsurance. In the following tables, the demand-change rate is the percentage of perturbations that altered the insurer's demand, and the improving-deviation rate is the percentage with \(\delta_{c,m}^{\,r}>0\). The median, upper percentiles, and maximum summarize the distribution of \(100\delta_{c,m}^{\,r}\), pooled across the reported solutions at each radius within each experimental setting.

\subsection{Results under the main experimental setting}

Tables~\ref{tab:supp-dist-main-no-cap} and~\ref{tab:supp-dist-main-cap} report the distributional summaries under the main experimental setting, without and with the affordability cap, respectively. No improving unilateral deviations were detected in either setting at any of the tested radii. At the smallest radius, \(r=0.1\), all sampled deviations altered the insurer's demand, yet the maximum observed relative objective changes remained negative: approximately \(-0.08\%\) without the affordability cap and \(-0.10\%\) with the cap. As a result, these tests examined nearby pricing deviations that produced actual changes in market participation. At \(r=5.0\), random sampling often generated large departures from the reported pricing vectors for which reinsurance reoptimization $\widetilde{\mcP}$ failed to identify feasible choices of \(A\) and \(M\). The resulting constraint-violation penalties in \(\bar\Phi\) explain the large negative relative objective changes reported in the table, particularly under the affordability cap.

\begin{table}[htp]
\centering
\caption{
Distributional summary of the relative objective changes
\(\delta_{c,m}^{\,r}\) without the affordability cap under the main
experimental setting. At each radius, 2,000 deviations per solution
are pooled across the ten equilibrium solutions reported in
Tables~3--4 of the main manuscript.
The median, percentiles, and maximum are computed for
\(100\delta_{c,m}^{\,r}\).
Demand-change and improving-deviation rates are also expressed
as percentages.
}
\label{tab:supp-dist-main-no-cap}
\resizebox{\textwidth}{!}{%
\begin{tabular}{c|rrrrrrr}
\hline
Radius
& Samples
& Demand change
& Improving deviations
& Median
& 95th pct.
& 99th pct.
& Maximum \\
\hline
\(0.1\)
& 20,000 & 100.00 & 0.00
& \(-0.65\) & \(-0.36\) & \(-0.27\) & \(-0.08\) \\
\(1.0\)
& 20,000 & 100.00 & 0.00
& \(-12.46\) & \(-4.88\) & \(-3.27\) & \(-1.11\) \\
\(5.0\)
& 20,000 & 100.00 & 0.00
& \(-6{,}268.03\) & \(-850.94\) & \(-268.18\) & \(-7.22\) \\
\hline
\end{tabular}%
}
\end{table}

\begin{table}[htp]
\centering
\caption{
Distributional summary of the relative objective changes
\(\delta_{c,m}^{\,r}\) with the affordability cap under the main
experimental setting, pooled across the ten equilibrium solutions
reported in Tables~5--6 of the main manuscript.
The sampling procedure, summary statistics, and units follow
Table~\ref{tab:supp-dist-main-no-cap}.
}
\label{tab:supp-dist-main-cap}
\resizebox{\textwidth}{!}{%
\begin{tabular}{c|rrrrrrr}
\hline
Radius
& Samples
& Demand change
& Improving deviations
& Median
& 95th pct.
& 99th pct.
& Maximum \\
\hline
\(0.1\)
& 20,000 & 100.00 & 0.00
& \(-3.45\) & \(-0.99\) & \(-0.61\) & \(-0.10\) \\
\(1.0\)
& 20,000 & 100.00 & 0.00
& \(-32{,}492.86\) & \(-99.21\) & \(-25.65\) & \(-6.91\) \\
\(5.0\)
& 20,000 & 100.00 & 0.00
& \(-581{,}004.37\) & \(-96{,}334.65\)
& \(-31{,}722.37\) & \(-61.73\) \\
\hline
\end{tabular}%
}
\end{table}

\subsection{Results with a larger acceptance threshold}

Tables~\ref{tab:supp-dist-loose-no-cap}
and~\ref{tab:supp-dist-loose-cap} report the corresponding results
for solutions obtained with \(\tau_{\mathrm{acc}}=10^{-4}\).
For each setting, 2,000 perturbations were generated per solution
at each radius, yielding 20,000 samples per radius across the
ten reported solutions.

\begin{table}[htp]
\centering
\caption{
Distributional summary of the relative objective changes
\(\delta_{c,m}^{\,r}\) without the affordability cap using
\(\tau_{\mathrm{acc}}=10^{-4}\), pooled across the ten equilibrium
solutions reported in Table~\ref{tab:supp-equil}.
The sampling procedure, summary statistics, and units follow
Table~\ref{tab:supp-dist-main-no-cap}.
}
\label{tab:supp-dist-loose-no-cap}
\resizebox{\textwidth}{!}{%
\begin{tabular}{c|rrrrrrr}
\hline
Radius
& Samples
& Demand change
& Improving deviations
& Median
& 95th pct.
& 99th pct.
& Maximum \\
\hline
\(0.1\)
& 20,000 & 100.00 & 0.00
& \(-0.60\) & \(-0.34\) & \(-0.24\) & \(-0.07\) \\
\(1.0\)
& 20,000 & 100.00 & 0.00
& \(-11.49\) & \(-4.68\) & \(-3.13\) & \(-1.10\) \\
\(5.0\)
& 20,000 & 100.00 & 0.00
& \(-5{,}545.32\) & \(-689.76\) & \(-241.38\) & \(-5.83\) \\
\hline
\end{tabular}%
}
\end{table}

\begin{table}[htp]
\centering
\caption{
Distributional summary of the relative objective changes
\(\delta_{c,m}^{\,r}\) with the affordability cap using
\(\tau_{\mathrm{acc}}=10^{-4}\), pooled across the ten equilibrium
solutions reported in Table~\ref{tab:supp-equilibria-afford}.
The sampling procedure, summary statistics, and units follow
Table~\ref{tab:supp-dist-main-no-cap}.
}
\label{tab:supp-dist-loose-cap}
\resizebox{\textwidth}{!}{%
\begin{tabular}{c|rrrrrrr}
\hline
Radius
& Samples
& Demand change
& Improving deviations
& Median
& 95th pct.
& 99th pct.
& Maximum \\
\hline
\(0.1\)
& 20,000 & 100.00 & 1.02
& \(-3.24\) & \(-0.88\) & \(0.02\) & \(2.96\) \\
\(1.0\)
& 20,000 & 100.00 & 0.00
& \(-32{,}386.13\) & \(-99.31\) & \(-24.75\) & \(-6.29\) \\
\(5.0\)
& 20,000 & 100.00 & 0.00
& \(-599{,}332.55\) & \(-102{,}940.48\)
& \(-31{,}733.12\) & \(-344.60\) \\
\hline
\end{tabular}%
}
\end{table}

Without the affordability cap, no improving unilateral deviations were detected at any tested radius.
With the affordability cap, improving deviations were identified for two of the ten solutions at \(r=0.1\), with maximum observed relative objective gains of approximately \(2.96\%\) and \(0.89\%\). Together, these accounted for 204 of the 20,000 perturbations at this radius, corresponding to an improving-deviation rate
of \(1.02\%\). No improving deviations were detected for the remaining eight solutions or at the larger tested radii.

\printbibliography[title={References}]





\end{document}

\endinput

%% file: mycommands.tex
\makeatletter
\let\save@mathaccent\mathaccent
\newcommand*\if@single[3]{%
  \setbox0\hbox{${\mathaccent"0362{#1}}^H$}%
  \setbox2\hbox{${\mathaccent"0362{\kern0pt#1}}^H$}%
  \ifdim\ht0=\ht2 #3\else #2\fi
  }
\newcommand*\rel@kern[1]{\kern#1\dimexpr\macc@kerna}
\newcommand*\widebar[1]{\@ifnextchar^{{\wide@bar{#1}{0}}}{\wide@bar{#1}{1}}}
\newcommand*\wide@bar[2]{\if@single{#1}{\wide@bar@{#1}{#2}{1}}{\wide@bar@{#1}{#2}{2}}}
\newcommand*\wide@bar@[3]{%
  \begingroup
  \def\mathaccent##1##2{%
    \let\mathaccent\save@mathaccent
    \if#32 \let\macc@nucleus\first@char \fi
    \setbox\z@\hbox{$\macc@style{\macc@nucleus}_{}$}%
    \setbox\tw@\hbox{$\macc@style{\macc@nucleus}{}_{}$}%
    \dimen@\wd\tw@
    \advance\dimen@-\wd\z@
    \divide\dimen@ 3
    \@tempdima\wd\tw@
    \advance\@tempdima-\scriptspace
    \divide\@tempdima 10
    \advance\dimen@-\@tempdima
    \ifdim\dimen@>\z@ \dimen@0pt\fi
    \rel@kern{0.6}\kern-\dimen@
    \if#31
      \overline{\rel@kern{-0.6}\kern\dimen@\macc@nucleus\rel@kern{0.4}\kern\dimen@}%
      \advance\dimen@0.4\dimexpr\macc@kerna
      \let\final@kern#2%
      \ifdim\dimen@<\z@ \let\final@kern1\fi
      \if\final@kern1 \kern-\dimen@\fi
    \else
      \overline{\rel@kern{-0.6}\kern\dimen@#1}%
    \fi
  }%
  \macc@depth\@ne
  \let\math@bgroup\@empty \let\math@egroup\macc@set@skewchar
  \mathsurround\z@ \frozen@everymath{\mathgroup\macc@group\relax}%
  \macc@set@skewchar\relax
  \let\mathaccentV\macc@nested@a
  \if#31
    \macc@nested@a\relax111{#1}%
  \else
    \def\gobble@till@marker##1\endmarker{}%
    \futurelet\first@char\gobble@till@marker#1\endmarker
    \ifcat\noexpand\first@char A\else
      \def\first@char{}%
    \fi
    \macc@nested@a\relax111{\first@char}%
  \fi
  \endgroup
}
\makeatother

\DeclareMathOperator*{\argmax}{\arg\!\max}

\renewcommand*{~}{\relax\ifmmode\sim\else\nobreakspace{}\fi}

\newcommand{\BFg}{\bm{g}}

\newcommand{\BFs}{\bm{s}}

\newcommand{\BFx}{\bm{x}}

\newcommand{\BFlambda}{\bm{\lambda}}

\newcommand{\sfH}{\mathsf{H}}

\newcommand{\mcC}{\mathcal{C}}

\newcommand{\mcG}{\mathcal{G}}

\newcommand{\mcP}{\mathcal{P}}

\newcommand{\mcS}{\mathcal{S}}